\documentclass[reqno,10pt]{amsart}
\usepackage{epsfig}
\usepackage[hypertex]{hyperref}
\usepackage{graphicx}

\newtheorem{lem}{Lemma}[section]
\newtheorem{pro}{Proposition}[section]
\newtheorem{thm}{Theorem}[section]

\newtheorem{rmk}{Remark}[section]

\numberwithin{equation}{section}

\title[Asymptotic distribution related to mildly-explosive AR(2)]
{Asymptotic distributions related to mildly-explosive second order
autoregressive models}

\date{July 20, 2014.}

\begin{document}

\maketitle

\author{
\begin{center}
Hui Jiang, \; Mingming Yu
%\footnote{ E-mail: huijiang@nuaa.edu.cn}, Yu Mingming\footnote{E-mail: mengyilianmeng@163.com}
\\
\vskip 0.2cm
{\it Department of Mathematics, Nanjing University of Aeronautics and Astronautics
\\ Nanjing 210016, P.R.China\\
\vskip 0.1cm

huijiang@nuaa.edu.cn, \; mengyilianmeng@163.com}

\vskip 0.4cm

Guangyu Yang %\footnote{ E-mail: guangyu@zzu.edu.cn}
\\
\vskip 0.2cm
{\it School of Mathematics and Statistics, Zhengzhou University
\\ Zhengzhou 450001, P.R.China\\
\vskip 0.1cm
guangyu@zzu.edu.cn}
\end{center}
}

\begin{abstract}
In this paper, we consider the normalized least squares estimator of the parameter in a mildly-explosive first-order autoregressive model with dependent errors which are modeled as a mildly-explosive AR(1) process. We prove that the estimator has a Cauchy limit law which provides a bridge between moderate deviation asymptotics and the earlier results on the local to unity and explosive autoregressive models.
In particular, the results can be applied to understand the near-integrated second order autoregressive processes. Simulation studies are also carried out to assess the performance of least squares estimation in finite samples.

\vskip10pt


\vskip10pt

\end{abstract}

\section{Introduction}\label{1}

\noindent There is a lot of econometric literature over the last three decades, which has focused on the issue of testing for the unit root hypothesis in econometric time series. Regression asymptotics with roots at or near unity have played an important role in time series econometrics. In order to cover more general time series structure, it has become popular in econometric methodology to study the models which permit that the regressors and the errors have substantial heterogeneity and dependence over time. In this paper, we mainly analyse a dynamic first order autoregressive model which the errors are dependent.
More precisely, we consider the following autoregressive model driven by a autoregressive error,
\begin{equation}\label{model}
\begin{aligned}
\begin{cases}
X_{k,n}={\theta_n}X_{k-1,n}+\varepsilon_{k,n}\\
\varepsilon_{k,n}={\rho_n}\varepsilon_{k-1,n}+ V_{k}
\end{cases}
\end{aligned},\qquad k=1,2,\dots,n,\;n\geq1,
\end{equation}
where the parameters $\theta_n$ and $\rho_n$ are unknown, $(X_{k,n})_{1\leq k\leq n}$ is observed, and the noise $(V_k)_{k\geq1}$ is a sequence of independent and identically distributed (i.i.d.) random variables with zero mean and a finite variance $\sigma^2$. For convenience, let $X_{0,n}=\varepsilon_{0,n}:=0$ for every $n$.
It is well-known that the least squares estimator of the parameter $\theta_n$ based on the observations $(X_{k,n})_{1\leq k\leq n}$ can be given by
\begin{align}
\hat{\theta}_n=\frac{\sum_{k=1}^n{X_{k,n} X_{k-1,n}}}{\sum_{k=1}^n{X_{k-1,n}^2}}.
\end{align}
To obtain the estimator of the parameter $\rho_n$, we can replace $\theta_n$ by $\hat{\theta}_n$ in (\ref{model}), and denote the estimators of the errors $(\varepsilon_{k,n})$ by
\begin{align}
\hat{\varepsilon}_{k,n}=X_{k,n}-\hat{\theta}_nX_{k-1,n},\qquad k=1,2,\ldots,n,
\end{align}
then the least squares estimators of $\rho_n$ can be defined as
\begin{align}
\hat{\rho}_n=\frac{\sum_{k=1}^n{\hat{\varepsilon}_{k,n} \hat{\varepsilon}_{k-1,n}}}{\sum_{k=1}^n{\hat{\varepsilon}_{k-1,n}^2}},
\end{align}
where $\hat{\varepsilon}_{0,n}:=0$ for every $n$.

\vskip5pt

The model (\ref{model}) has a close connection with some existing models. Firstly, we fix the autoregressive coefficient $\theta_n$, i.e. let $\theta_n\equiv\theta$. If $\rho_n\equiv0$, then the model (\ref{model}) turns to be the classic autoregressive process with i.i.d. errors. In this case, the asymptotic behaviors of $\hat{\theta}_n$ have been examined thoroughly. For example, when the model is {\it stationary} ($|\theta|<1$), under some moment conditions, Anderson \cite{Anderson}
showed asymptotic normality of $\hat{\theta}_n-\theta$.
However, as pointed out previously by Anderson \cite{Anderson}, White \cite{White}, and Dickey \& Fuller \cite{DickFull1979}, the situation becomes more complicated for the {\it critical} case ($|\theta|=1$) and the {\it explosive} case ($|\theta|>1$), where the limiting distributions are functionals of Brownian motion and standard Cauchy, respectively. In addition, if the regressive coefficient $\rho_n$ in the errors is also fixed, i.e. $\rho_n\equiv\rho$, to answer some open problems on the Durbin-Watson statistic, Bercu and Pro\"{i}a \cite{Bercu-2013} investigated the asymptotic normality of the least squares estimators $\hat{\theta}_n$ and $\hat{\rho}_n$, while Bitseki Penda {\it et al.} \cite{Penda} studied the moderate deviations, both in the {\it stationary} cases, i.e. $|\theta|<1$ and $|\rho|<1$.

\vskip5pt

Note that, the above investigations can capture the
phenomena with phase transition type characteristics, i.e.
from the {\it stationary} to the {\it critical}, and from the {\it critical}
to the {\it explosive}, which just corresponds to the transition of the limiting distribution, from the {\it normal} distribution to the {\it functional of Brownian motion}, to the standard {\it Cauchy} distribution.

\vskip5pt

To understand this phase transition and handle the data that allows for large shocks
in the dynamic structure of the model, more recently, some attention has been
dedicated to the autoregressive models with the dynamic coefficient.
To accommodate this observation, $\theta_n$ is allowed to depend on the sample size $n$. Similar to the above mentioned, first let the noise $(\varepsilon_{k,n})_{1\leq k\leq n}$ in (\ref{model}) be a sequence of i.i.d. random variables, i.e. $\rho_n\equiv0$. Recall that, we say that the model (\ref{model}) is local to unity if $\theta_n=1+c/n$. This model has been proved useful in analysing the near integrated processes, in establishing the local asymptotic properties of tests and in the construction of confidence intervals.
Chan \& Wei \cite{Chan-Wei} and Phillips \cite{Phillips} showed that the asymptotic distribution of $\hat{\theta}_n-\theta_n$ is some kind of functional of Brownian motion.
To characterize great deviations from unity and understand the phase transition, Phillips \& Magdalinos \cite{Philllips-Magdalinos} considered the case, $\theta_n=1+c/\kappa_n$, where $(\kappa_n)_{n\geq1}$ is a deterministic sequence increasing to infinity satisfying $\kappa_n=o(n)$ and it represents moderate deviations from unity. They showed that, $\hat{\theta}_n-\theta_n$ has a $\sqrt{n\kappa_n}$ rate of convergence and a asymptotic normal distribution when $c<0$, and $\hat{\theta}_n-\theta_n$ has a $\kappa_n\theta_n^n$ rate of convergence and a standard Cauchy limit distribution when $c>0$. More interestingly, their results {\it match} the standard limit theory of the {\it fixed} coefficients model and partially bridge the {\it stationary}, the {\it local to unity} and the {\it explosive} cases. Very recently, Miao {\it et al.} \cite{Miao-Yang} derived the moderate deviations of $\hat{\theta}_n-\theta_n$ as $\theta_n\to1$ within stationary regions, which also {\it matches} the standard limit theory of the {\it fixed} coefficient model. While, if the noise $(\varepsilon_{k,n})_{1\leq k\leq n}$ in (\ref{model}) has dependent structure, one can refer to Giraitis \& Phillips \cite{Giraitis} for the martingale difference noise, Phillips \& Magdalinos \cite{Philllips-Magdalinos-1}, and Magdalinos \cite{Magdalinos2012} for
some weakly and strongly dependent noises.

\vskip5pt

It is remarkable that, to provide a better asymptotic framework for the nearly integrated
first order autoregressive model driven by an AR(1) process with root approaching the unity, Nabeya \& Perron \cite{Nabeya} also introduced the model (\ref{model}), where
they put $\theta_n=1+\gamma_1/n$ and $\rho_n=1+\gamma_2/n$. And they showed that the asymptotic distribution of $\hat{\theta}_n-\theta_n$ is some kind of functional of Brownian motion. In fact, just as pointed out by Nabeya \& Perron \cite{Nabeya}, the model (\ref{model}) can also be regarded as an approximate version of the second order autoregressive process with two unit roots. For more detailed explanations on this model, please refer to Nabeya \& Perron \cite{Nabeya}, Hasza and Fuller \cite{HaszaFuller1979}, or Chan \cite{Chan2009}.

\vskip5pt

Then, motivated by the above discussions, we will devote to the asymptotic properties of $\hat{\theta}_n-\theta_n$ in the nearly integrated first order autoregressive model driven by the nearly integrated AR(1) process. In the present paper, we mainly consider the case, $|\theta_n|\to1$ and $|\rho_n|\to1$ both within the {\it explosive} regions. To be specific, when $|\theta_n|=1+\gamma_1/k_n$ and $|\rho_n|=1+\gamma_2/k_n$, where $\gamma_1,\gamma_2>0$ and $(k_n)_{n\geq1}$ is a sequence of positive numbers increasing to infinity at a rate slower than $n$, we prove that the limiting distribution of the least squares estimator $\hat{\theta}_n-\theta_n$ is the Cauchy distribution which partially matches the standard limit theory of the aforementioned models. Just as pointed out previously by Phillips \& Magdalinos \cite{Philllips-Magdalinos-1}, and Magdalinos \cite{Magdalinos2012}, the resulting Cauchy limit law for the normalized autoregressive coefficient suggests that the limit theory is invariant to the dependence structure of the innovation errors in the {\it explosive} case. However, there also appear some interesting phenomena when $\theta_n$ and $\rho_n$ have the different signs. In the other preprint \cite{JYY}, we mainly analysed the case, $|\theta_n|\to1$ and $|\rho_n|\to1$ both within the {\it stationary} regions, and obtained the asymptotic normality and moderate deviations of the least squares estimators $\hat{\theta}_n$, $\hat{\rho}_n$ and the Durbin-Watson statistic. Finally, it is worthwhile to note that our results can be applied to understand the near-integrated second order autoregressive process.

\vskip5pt

The rest of this paper is organized as follows. The next section is devoted to the descriptions of our main results and some related discussions. In Section \ref{3}, we carry out some statistical simulations for the main results which imply
that our asymptotic results well match the
finite-sample properties of the estimators. Then, the technical proofs of main results are completed in the remaining sections.

\section{Results and discussions}\label{2}

\subsection{Main results} The following are our main results.

\begin{thm}\label{thm-1}
For model (\ref{model}) with $\theta_n=1+\gamma_1/k_n$, $\rho_n=1+\gamma_2/k_n$, and $k_n=o(n)$, we have, as $n\to\infty$,
\vskip3pt

\noindent(1) if $\gamma_1>\gamma_2>0$, then
$$
\frac{\gamma_1+\gamma_2}{2\gamma_1(\gamma_1-\gamma_2)}
k_n\theta_n^n\rho_n^{-n}(\hat\theta_n-\theta_n)\stackrel{\mathcal{L}}{\longrightarrow}\frac {\xi_{\rho}}{\xi_{\theta}};
$$

\noindent(2) if $\gamma_2>\gamma_1>0$, then
$$
\frac{\gamma_1+\gamma_2}{2(\gamma_2-\gamma_1)}
\rho_n^n \theta_n^{-n}\left(k_n(\hat\theta_n-\theta_n)-(\gamma_2-\gamma_1)\right)
\stackrel{\mathcal{L}}{\longrightarrow}\frac {\xi_{\theta}}{\xi_{\rho}},
$$
where $\stackrel{\mathcal{L}}\longrightarrow$ denotes the convergence in distribution
and $(\xi_{\theta}, \xi_{\rho})\sim N(0, \Lambda)$ with the covariance matrix
$$
\Lambda=
\begin{pmatrix}
\frac{\sigma^2}{2\gamma_1} & \frac{\sigma^2}{\gamma_1+\gamma_2}
\\
\frac{\sigma^2}{\gamma_1+\gamma_2} & \frac{\sigma^2}{2\gamma_2}
\end{pmatrix};
$$

\noindent(3) if $\gamma_1=\gamma_2=\gamma>0$, then
$$
\frac{n}{k_n}\left(n(\hat\theta_n-\theta_n)-\theta_n\right)
\stackrel{\mathcal{L}}{\longrightarrow}\frac {\zeta_{\theta}}{\varphi_{\theta}},
$$
where $\left(\varphi_{\theta},\zeta_{\theta}\right)\sim N(0,\Xi)$ with the covariance matrix
$$
\Xi=
\begin{pmatrix}
\frac{\sigma^2}{2\gamma} & \frac{\sigma^2}{2\gamma^2}
\\
\frac{\sigma^2}{2\gamma^2}  & \frac{5\sigma^2}{8\gamma^3}
\end{pmatrix}.
$$
\end{thm}

Note that the asymptotic distributions in Theorem \ref{thm-1} are Cauchy distributions {\it centered} at $2\gamma_1/(\gamma_1+\gamma_2)$, $2\gamma_2/(\gamma_1+\gamma_2)$ and $1/\gamma$, respectively. However, it is surprised that the asymptotic distribution is a {\it standard} Cauchy distribution
when the parameters $\theta_n$ and $\rho_n$ have different signs.

\begin{thm}\label{thm-2}
For model (\ref{model}) with $\theta_n=1+\gamma_1/k_n$, $\rho_n=-1-\gamma_2/k_n$, and $k_n=n^\alpha$ for $\alpha\in(0,1)$, we have, as $n\to\infty$,
\vskip3pt

\noindent(1) if $\gamma_1>\gamma_2>0$, then
\[
\frac{1}{2\gamma_1}\sqrt{\frac{\gamma_2}{\gamma_1}}
k_n\theta_n^n\rho_n^{-n}(\hat\theta_n-\theta_n)\stackrel{\mathcal{L}}{\longrightarrow}C;
\]

\noindent(2) if $\gamma_2>\gamma_1>0$, then
\[
\frac{1}{2\gamma_2}\sqrt{\frac{\gamma_1}{\gamma_2}}
k_n\rho_n^n\theta_n^{-n}(\hat\theta_n-\rho_n)
\stackrel{\mathcal{L}}{\longrightarrow}C,
\]
where $C$ denotes the standard Cauchy random variable.
%\[
%\Upsilon=\left(
%  \begin{array}{cc}
%    \frac{\sigma^2}{2\gamma_1} & 0 \\
%    0 & \frac{\sigma^2}{2\gamma_2} \\
%  \end{array}
%\right).
%\]
\end{thm}

\begin{rmk}
When $\gamma_1>\gamma_2>0$, i.e. $|\theta_n|>|\rho_n|$, the least squares estimator $\hat{\theta}_n$ is the consistent estimator of $\theta_n$ both in Theorems \ref{thm-1} and \ref{thm-2}. However, it is mysterious that, when $\gamma_2\geq\gamma_1>0$, i.e. $|\theta_n|\leq|\rho_n|$, $\hat{\theta}_n$ has an asymptotic bias which is similar to the results in Bercu \& Pro\"{i}a \cite{Bercu-2013}, Stocker \cite{Stocker}, and Phillips \& Magdalinos \cite{Philllips-Magdalinos-1} in the stationary and near-stationary cases. In fact, we have
\begin{align}
k_n(\hat{\theta}_n-\theta_n)\stackrel{P}{\longrightarrow}\gamma_2-\gamma_1, \qquad if \; \rho_n >\theta_n>1,
\end{align}
and
\begin{align}
\hat{\theta}_n-\theta_n\stackrel{P}{\longrightarrow}-2, \qquad if \; \rho_n<-\theta_n<-1,
\end{align}
where $\stackrel{{P}}\longrightarrow$ denotes the convergence in probability. Please refer to Proposition \ref{lem-5} and Appendix for more details.
\end{rmk}

\begin{rmk}\label{link1}
For model (\ref{model}) with $\theta_n=1+c/n^\alpha$ for some $c>0$ and $\alpha\in(0,1)$,
Phillips \& Magdalinos \cite{Philllips-Magdalinos-1} considered some weakly dependent errors, i.e. $\varepsilon_{k,n}=\sum_{j=0}^{\infty}c_jV_{k-j}$, where the non-random sequence $(c_j)_{j\geq0}$ is independent of $n$.
Under some summability conditions on $(c_j)_{j\geq0}$, the asymptotic distribution of $\hat{\theta}_n-\theta_n$ is proved to be standard Cauchy.
However,  because $\xi_{\theta}$ and $\xi_{\rho}$  are not independent as well as $\zeta_{\theta}$ and $\varphi_{\theta}$, it is interesting that the limiting distributions of $\hat{\theta}_n-\theta_n$ are not standard Cauchy distribution as shown in our Theorem \ref{thm-1}. As mentioned earlier, Theorem \ref{thm-2} shows that the limiting distribution of the normalization of $\hat{\theta}_n$ turns to be standard Cauchy which matches the results in White \cite{White}, Anderson \cite{Anderson}, Phillips \& Magdalinos \cite{Philllips-Magdalinos}, \cite{Philllips-Magdalinos-1}, and Magdalinos \cite{Magdalinos2012}.
Statistical simulations in Section \ref{3} also illustrate these.
\end{rmk}

\subsection{Discussions}\label{2.2}

It is still worthwhile to give some additional comments on our results and other related problems.

\begin{enumerate}
  \item In fact, under an additional symmetry assumption on the distribution of the noise $(V_k)_{k\geq1}$, Theorem \ref{thm-1} still holds in the case, $\theta_n\to-1$ and $\rho_n\to-1$, both within the {\it explosive} regions. Suppose that
  \begin{align}
  \left\{
  \begin{array}{ll}
  X_{k,n}&={\theta_n}X_{k-1,n}+\varepsilon_{k,n},\\
  \varepsilon_{k,n}&={\rho_n}\varepsilon_{k-1,n}+V_k
  \end{array},
  \right.\quad k=1,2,\ldots,n,\;n\geq1,
  \end{align}
  where the unknown parameters
  \[
  \theta_n=-1-\frac{\gamma_1}{k_n},\quad \rho_n=-1-\frac{\gamma_2}{k_n},\qquad
  \gamma_1>0,\;\gamma_2>0,
  \]
  and $(V_k)_{k\geq1}$ is a sequence of i.i.d. random variables with a {\it symmetric} distribution. Denote
  \begin{align*}
  &\alpha_n=-\theta_n,\quad \beta_n=-\rho_n,\\
  Y_{k,n}=(-1)^{k}X_{k,n},&\quad \eta_{k,n}=(-1)^k\varepsilon_{k,n},\quad W_k=(-1)^kV_k,
  \end{align*}
  then $(W_k)_{k\geq1}$ is a sequence of i.i.d. random variables with the same common distribution as that of $V_1$, and
  \begin{align}
  \left\{
  \begin{array}{ll}
  Y_{k,n}&={\alpha_n}Y_{k-1,n}+\eta_{k,n},\\
  \eta_{k,n}&={\beta_n}\eta_{k-1,n}+W_k
  \end{array},
  \right.\quad k=1,2,\ldots,n,\;n\geq1.
  \end{align}
  Putting
  \begin{align*}
  \hat{\theta}_n&=\frac{\sum_{k=1}^n{X_{k,n} X_{k-1,n}}}{\sum_{k=1}^n{X_{k-1,n}^2}},\quad
  \hat{\alpha}_n=\frac{\sum_{k=1}^n{Y_{k,n} Y_{k-1,n}}}{\sum_{k=1}^n{Y_{k-1,n}^2}},
  \end{align*}
  it is easy to see that $\hat{\alpha}_n=-\hat{\theta}_n$,
  hence, if the corresponding assumptions are satisfied, then Theorem \ref{thm-1}
  holds for the least squares estimator $\hat{\alpha}_n$, hence also for $\hat{\theta}_n$, except for some minor changes, i.e. the Cauchy limit distributions are centered at $-2\gamma_1/(\gamma_1+\gamma_2)$, $-2\gamma_2/(\gamma_1+\gamma_2)$ and $-1/\gamma$ respectively, and the
  removed term turns to be $\gamma_1-\gamma_2$ in the case of $\gamma_2>\gamma_1>0$.
  As for the other case,
  \[
  \theta_n=-1-\frac{\gamma_1}{k_n}, \quad \rho_n=1+\frac{\gamma_2}{k_n},\qquad \gamma_1>0, \; \gamma_2>0,
  \]
  by the same method, we can show that Theorem \ref{thm-2} also holds if the corresponding assumptions are satisfied.

  \item

  As mentioned in Remark \ref{link1}, Theorems \ref{thm-1} and \ref{thm-2} relate to the earlier work (White \cite{White}, Anderson \cite{Anderson}, Basawa \& Brockwell \cite{BasBro1984}, Nabeya \& Perron \cite{Nabeya}, Phillips \& Magdalinos \cite{Philllips-Magdalinos}, \cite{Philllips-Magdalinos-1}, Magdalinos \cite{Magdalinos2012}) on the {\it explosive} AR(1) process. For the Gaussian first order autoregressive model with fixed coefficient $|\theta|>1$, White proved that
  \begin{align}\label{White}
  \frac{\theta^n}{\theta^2-1}(\hat{\theta}_n-\theta)\stackrel{\mathcal{L}}
  {\longrightarrow}C,
  \end{align}
  where $C$ denotes the {\it standard} Cauchy random variable. Phillips \& Magdalinos showed that (\ref{White}) still holds, if the parameter $\theta$ and the Gaussian errors are respectively replaced by $\theta_n=1+\gamma/n^{\alpha}$\,($\gamma>0,\alpha\in(0,1)$), and i.i.d. (even some long range dependent) errors.
  However, Theorems \ref{thm-1} and \ref{thm-2} say that it also can be extended to some strongly dependent cases. This provides further evidence that the asymptotic theory is invariant to the dependence structure of the innovation errors in the {\it explosive} case.

  \item
  Note that $|\theta_n|\to1$ and $|\rho_n|\to1$, both within the explosive regions, hence our main results, Theorems \ref{thm-1} and \ref{thm-2}, maybe provide a bridge between those for unit root (or local to unity) processes and those that under the  explosive case with strongly dependent errors. Assume that $\gamma_1,\gamma_2>0$ and $k_n=n^\alpha$ for some $\alpha\in(0,1)$. Parts (1) of Theorems \ref{thm-1} and \ref{thm-2} become
  \begin{align}
  \frac{\gamma_1+\gamma_2}{2\gamma_1(\gamma_1-\gamma_2)}
  n^\alpha\theta_n^n\rho_n^{-n}(\hat\theta_n-\theta_n)
  \stackrel{\mathcal{L}}{\longrightarrow}C_1,
  \end{align}
  and
  \begin{align}
  \frac{1}{2\gamma_1}\sqrt{\frac{\gamma_2}{\gamma_1}}
  n^\alpha\theta_n^n\rho_n^{-n}(\hat\theta_n-\theta_n)
  \stackrel{\mathcal{L}}{\longrightarrow}C,
  \end{align}
  where $C_1$ denotes the Cauchy variate centered at $2\gamma_1/(\gamma_1+\gamma_2)$.
  It is notable that, ignoring the multiplicative constants, the convergence rate takes values in $(n,(\frac{1+\gamma_1}{1+\gamma_2})^n)$ as $\alpha$ ranges form $1$ to $0$. When $\alpha=0$, the model (\ref{model}) becomes a standard second order autoregressive model with two explosive characteristic roots, $1+\gamma_1$ and $1+\gamma_2$, which had been considered by Anderson \cite{Anderson}.
  Thus, the convergence rate of the serial correlation coefficient covers the interval $(n,(\frac{1+\gamma_1}{1+\gamma_2})^n)$, establishing a link between the asymptotic behavior of local to unity and explosive autoregressive models. However, when $\alpha=1$, this is replaced by the following local to unity
  limit theory developed by Nabeya \& Perron \cite{Nabeya},
  \begin{align}\label{Nabeya-Perron}
  n(\hat{\theta}_n-\theta_n)\stackrel{\mathcal{L}}{\longrightarrow}
  \frac{Q_{\gamma_1}(J_{\gamma_2}(1))^2}{2\int_0^1Q_{\gamma_1}(J_{\gamma_2}(s))^2ds}
  -\gamma_1,
  \end{align}
  where $(B(t))$ is the standard Brownian motion, $J_{\gamma_2}(t)=\int_0^te^{\gamma_2(t-x)}dB(s)$ is an Ornstein-Uhlenbeck process, and $Q_{\gamma_1}(J_{\gamma_2}(t))$ is the weighted integral of the process $(J_{\gamma_2}(t))$,
  \begin{align}\label{Nabeya-Perron1}
  Q_{\gamma_1}(J_{\gamma_2}(t)):=\int_0^te^{\gamma_1(t-s)}J_{\gamma_2}(s)ds.
  \end{align}

  \item
  More meaningly, when $\gamma_1=\gamma_2=\gamma>0$, model (\ref{model}) is just the second order autoregression with common {\it near-explosive} roots, $1+\gamma/k_n$.
  Part (3) of Theorem \ref{thm-1} when $k_n=n^\alpha$ for some $\alpha\in(0,1)$, turns to be
  \begin{align}
  n^{1-\alpha}\left(n(\hat\theta_n-\theta_n)-\theta_n\right)
  \stackrel{\mathcal{L}}{\longrightarrow}C_2,
  \end{align}
  where $C_2$ is a Cauchy variate centered at $1/\gamma$. The convergence rate covers a more smaller interval $(n,n^2)$ as $\alpha$ ranges form $1$ to $0$. When $\alpha=1$, it is natural to consider the local to unity limit theory (\ref{Nabeya-Perron}), however, when $\alpha=0$, Nielsen \cite{Nielsen2009} showed that the least squares estimator is inconsistent. Phillips \& Magdalinos
  \cite{Philllips-Magdalinos-2} provided a co-explosive system extension and an illustrative examples to explain the finding. And they also gave a consistent instrumental variable procedure. In addition, they pointed out that the least squares estimator is again consistent when $\theta_n\to1$ within the explosive region.

  \item
  As mentioned before, model (\ref{model}) can be regarded as a second order autoregressive process with two characteristic roots, $1+\gamma_1/k_n$ and $1+\gamma_2/k_n$. The present paper and the preprint \cite{JYY} systematically study the case, $\gamma_1\gamma_2>0$.
  It is natural to ask what will happen if $\gamma_1\gamma_2\leq0$, which had been studied by Rao \cite{Rao1961} when $k_n\equiv1$ and a root exceeding one and the other less than one in absolute value. More generally, if $\theta_n\to1$ and $\rho_n\to1$ with different rates, can we say something? To our knowledge, Phillips \& Lee \cite{PhiLee2012} recently have developed some limit theory for the nonstationary vector autoregression with mixed roots in the vicinity of unity.
\end{enumerate}

\section{Simulation studies}\label{3}

\noindent To further illustrate our main results, Theorems \ref{thm-1} and \ref{thm-2}, and understand the discussions in Section \ref{2.2}, using R software with the help of Jianbin Zhao, in this section we carry out some statistical simulations to examine the performance of the asymptotic results in finite samples. The results show that the limiting distributions match well with the finite samples distributions and the limiting distributions of $\hat{\theta}_n$, given that $\theta_n=1+\gamma_1/k_n,\;\rho_n=1+\gamma_2/k_n$, or,  $\theta_n=1+\gamma_1/k_n,\;\rho_n=-1-\gamma_2/k_n$, are respectively identical and equal to the mirror images of the limiting distributions of $\hat{\theta}_n$ given that $\theta_n=-1-\gamma_1/k_n$, $\rho_n=-1-\gamma_2/k_n$, or, $\theta_n=-1-\gamma_1/k_n,\;\rho_n=1+\gamma_2/k_n$, provided $\gamma_1,\;\gamma_2>0$. In addition, they also show that the Cauchy limit distributions are respectively biased and unbiased in Theorems \ref{thm-1} and \ref{thm-2}.

\vskip5pt

We now give some explanations for the simulations. Data are generated through model (\ref{model}) under the assumptions of Theorems \ref{thm-1} and \ref{thm-2}, where we let the noise $(V_k)_{k\geq1}$ be a sequence
of i.i.d. Gaussian random variables with zero
mean and unit variance. The sample size is $n=400$ and the number of replications is $1000$. In addition, we put $k_n=n^{1/3}$. In the following figures, the {\it blue} and {\it red curves} denote the density curves of Cauchy and finite samples distributions respectively. The first three groups correspond to parts (1), (2) and (3) of Theorem \ref{thm-1}. And the last two groups correspond to parts (1) and (2) of Theorem \ref{thm-2}.

\hspace*{5mm}\includegraphics[width=2in]{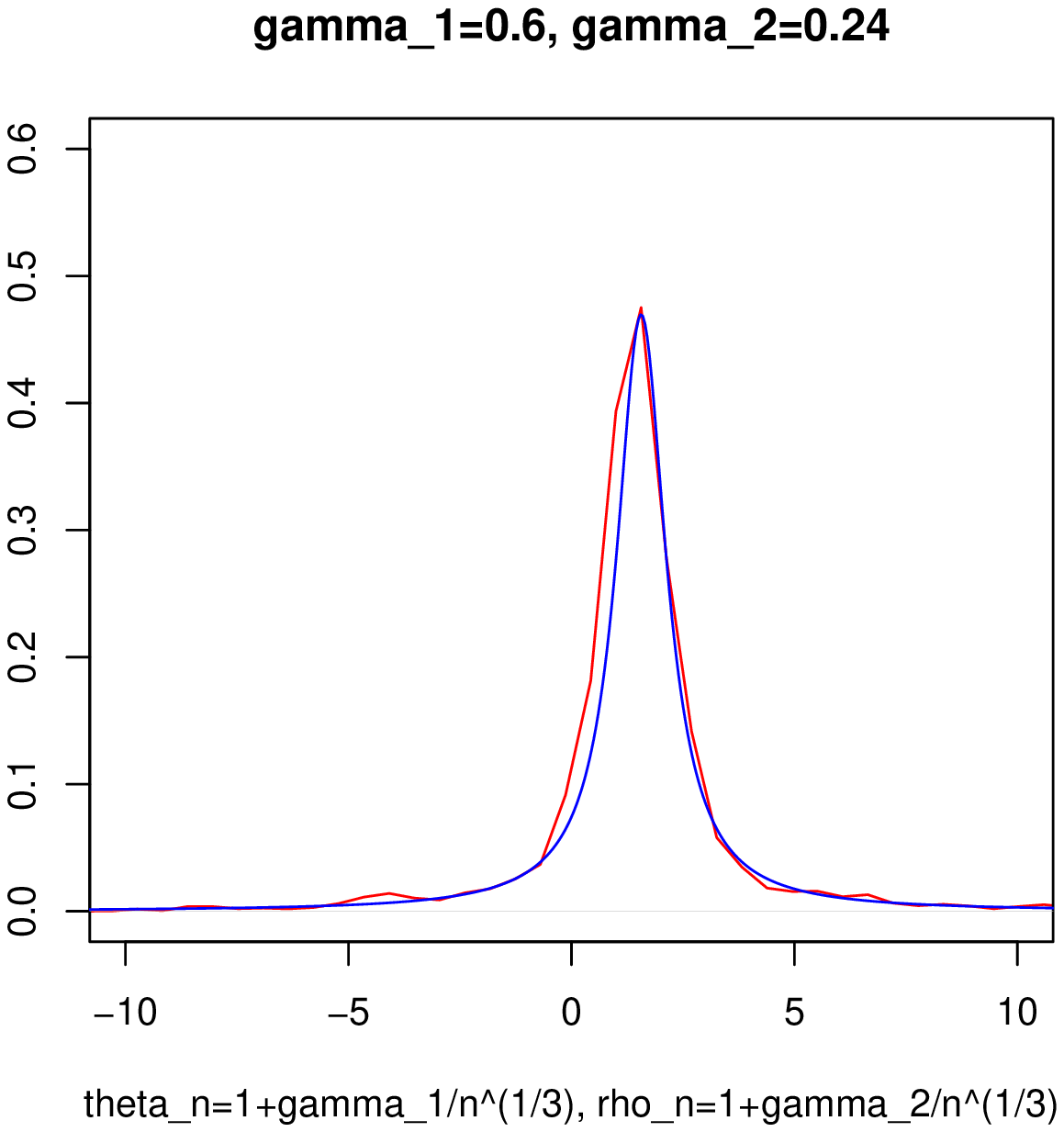}
\includegraphics[width=2in]{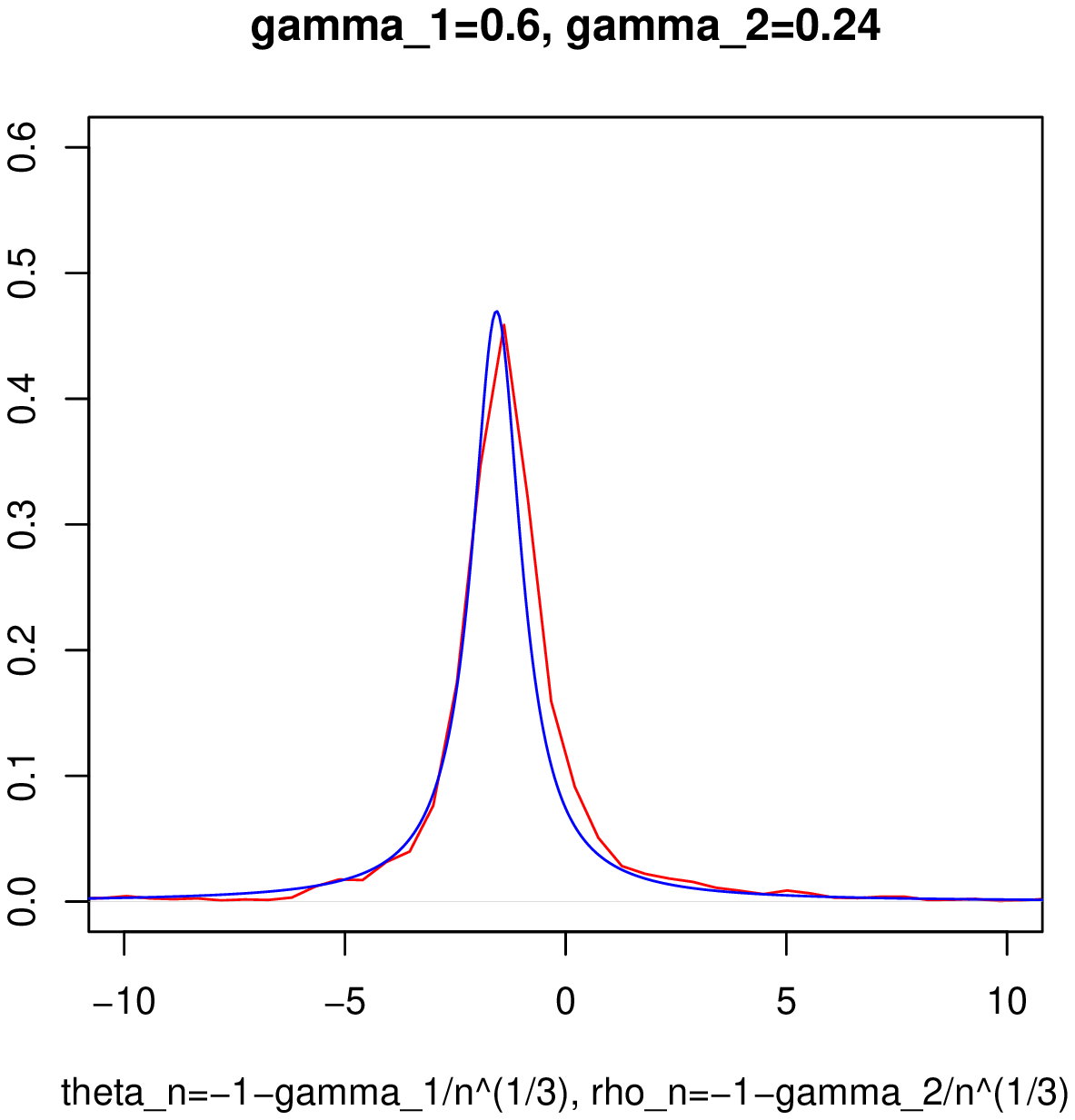}

\hspace*{5mm}\includegraphics[width=2in]{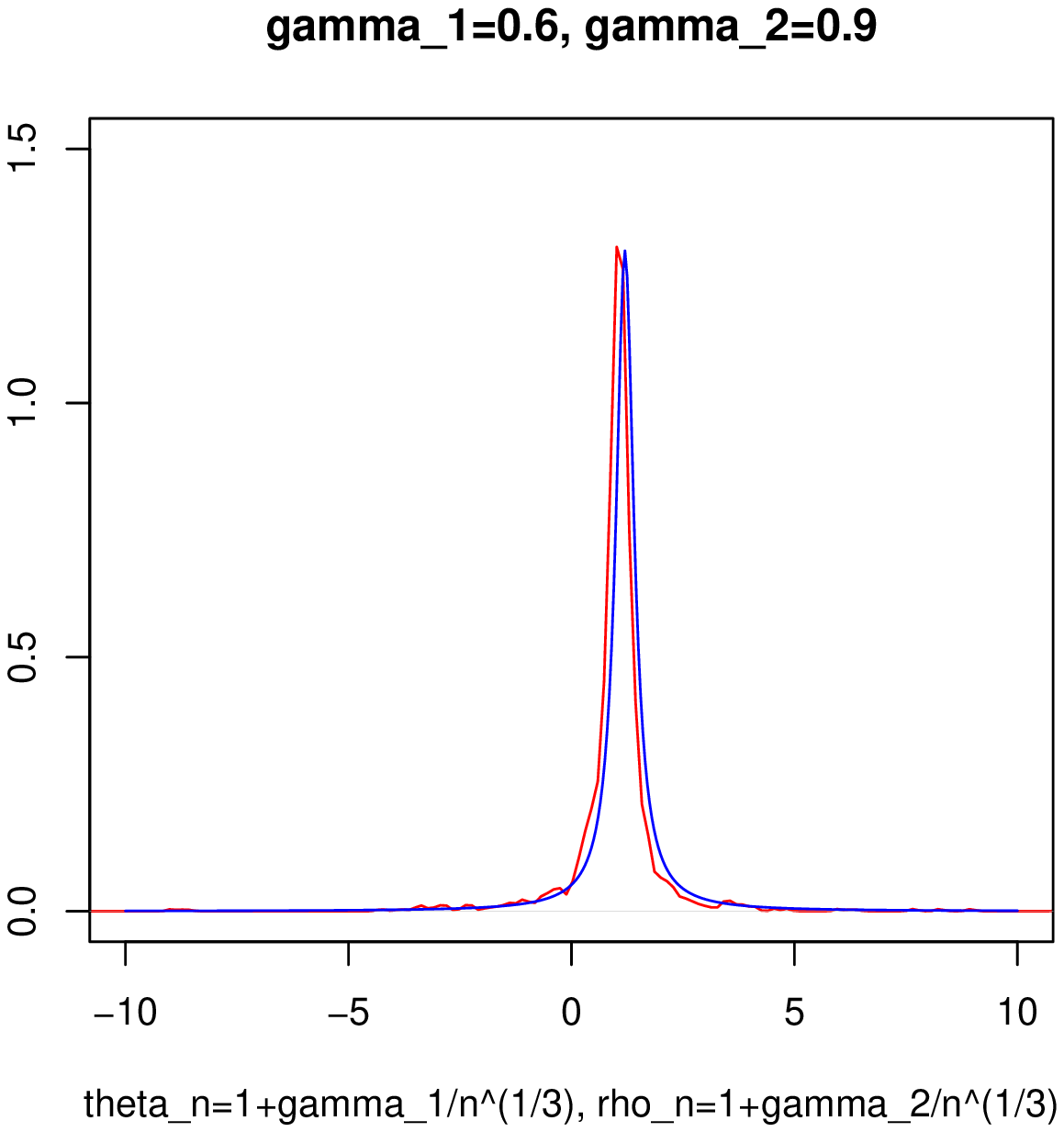}
\includegraphics[width=2in]{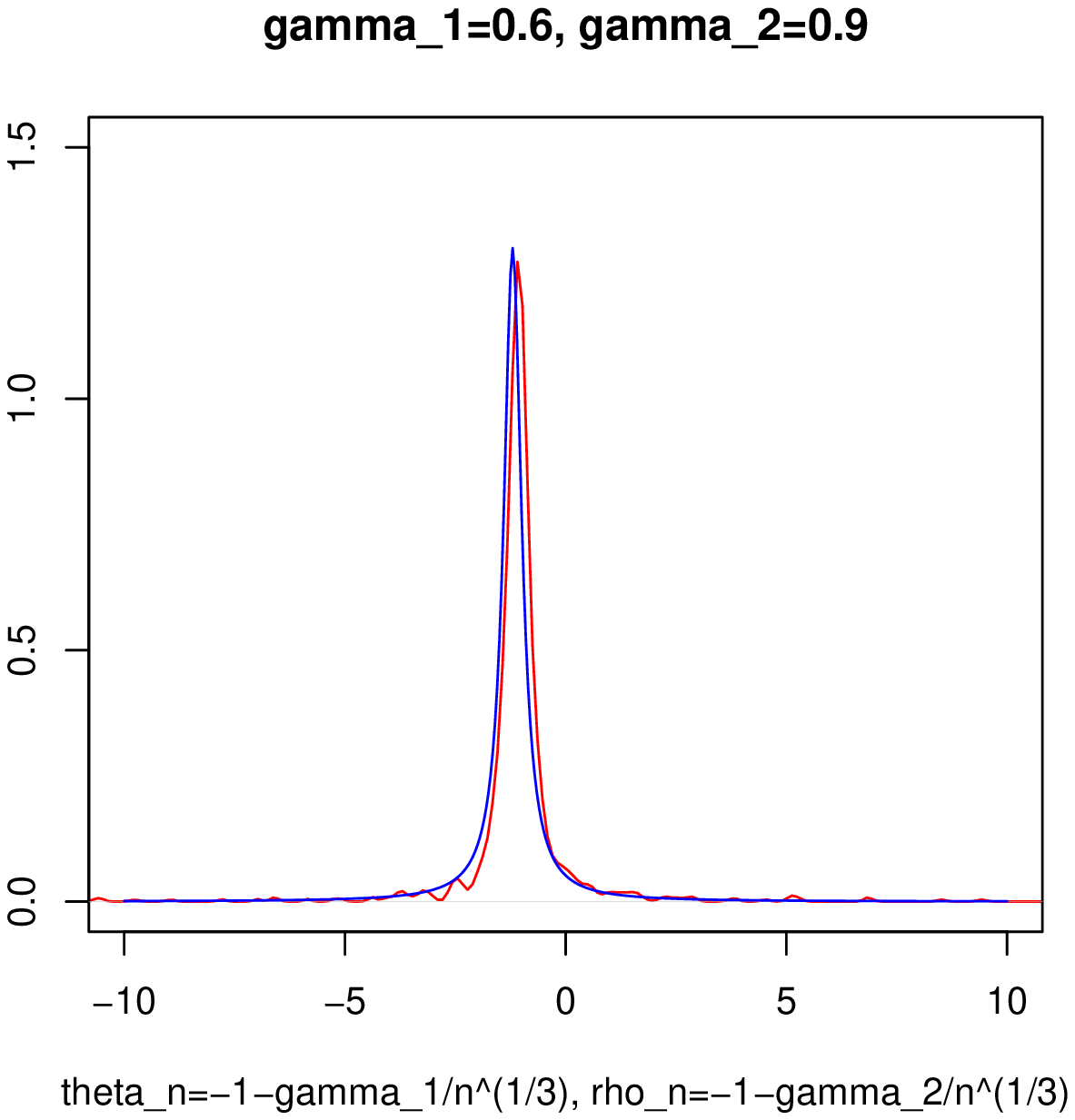}

\hspace*{5mm}\includegraphics[width=2in]{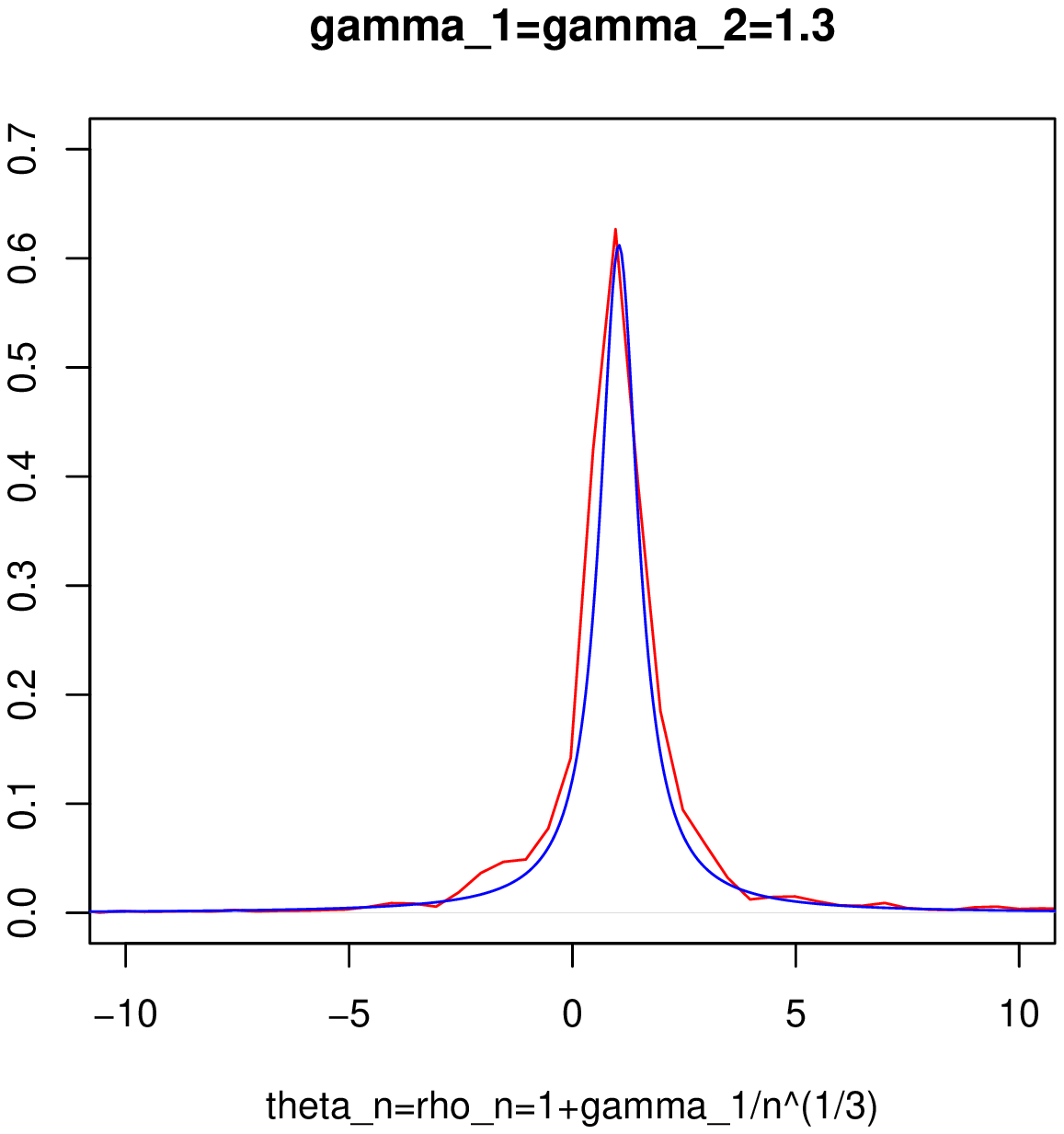}
\includegraphics[width=2in]{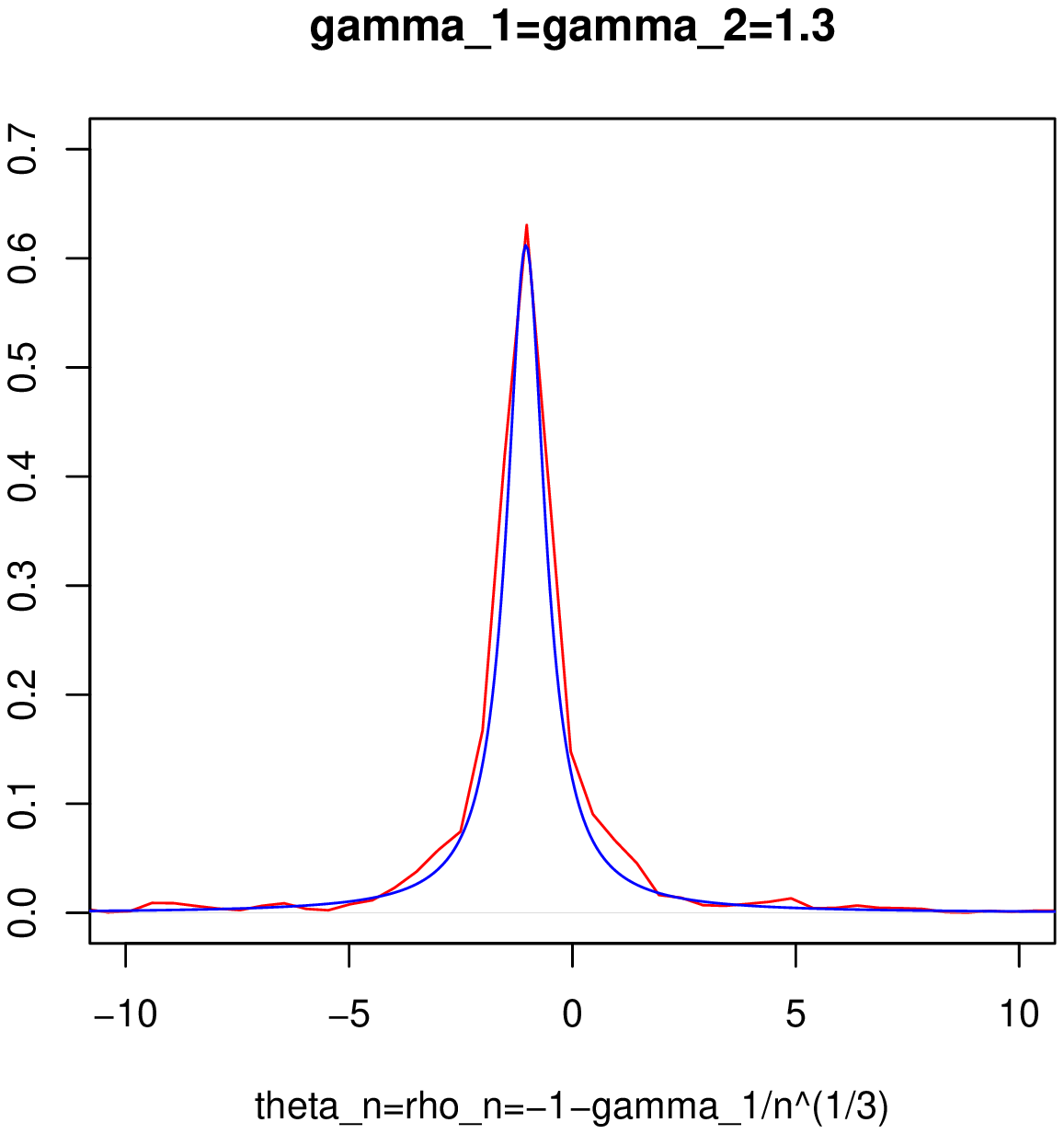}

\hspace*{5mm}\includegraphics[width=2in]{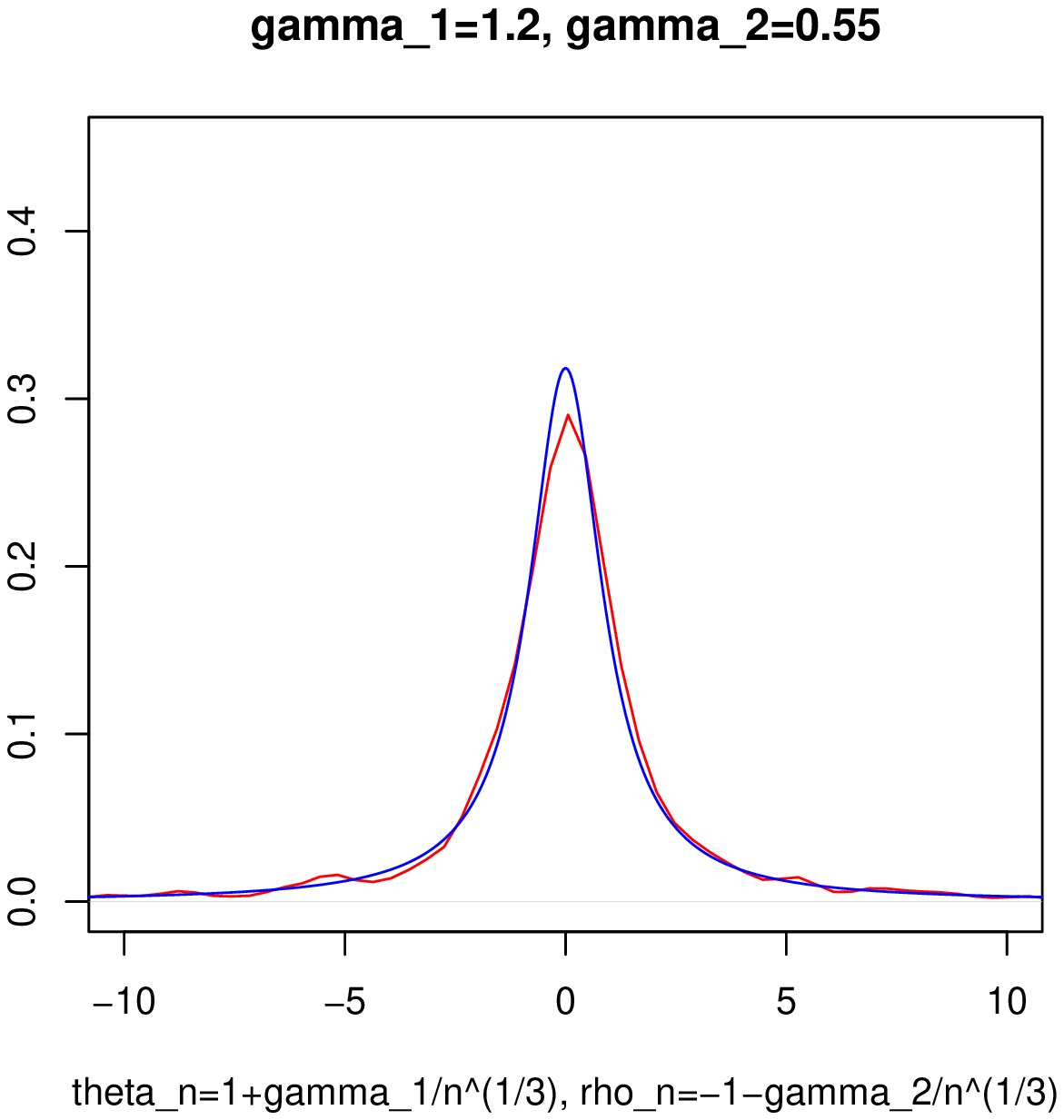}
\includegraphics[width=2in]{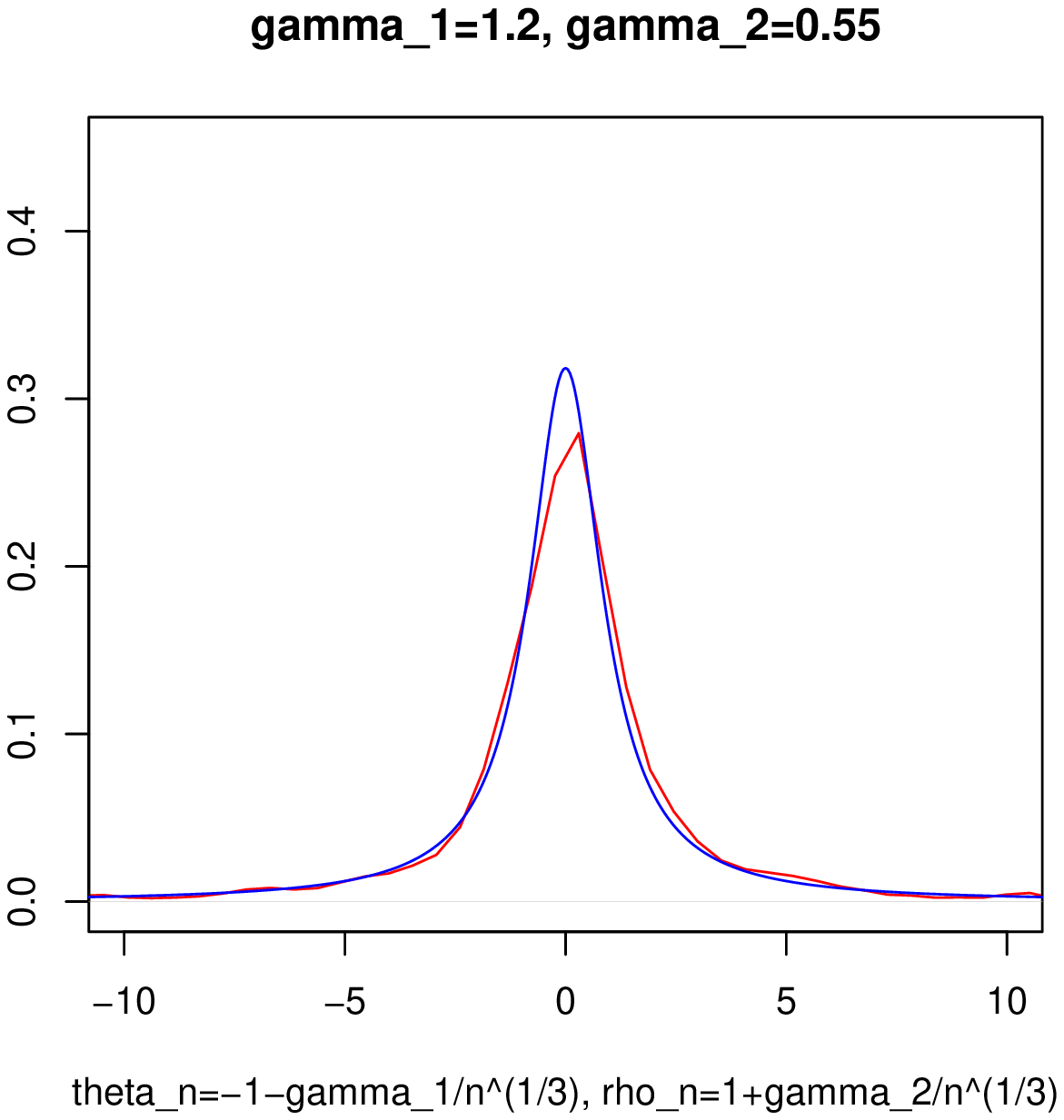}

\hspace*{5mm}\includegraphics[width=2in]{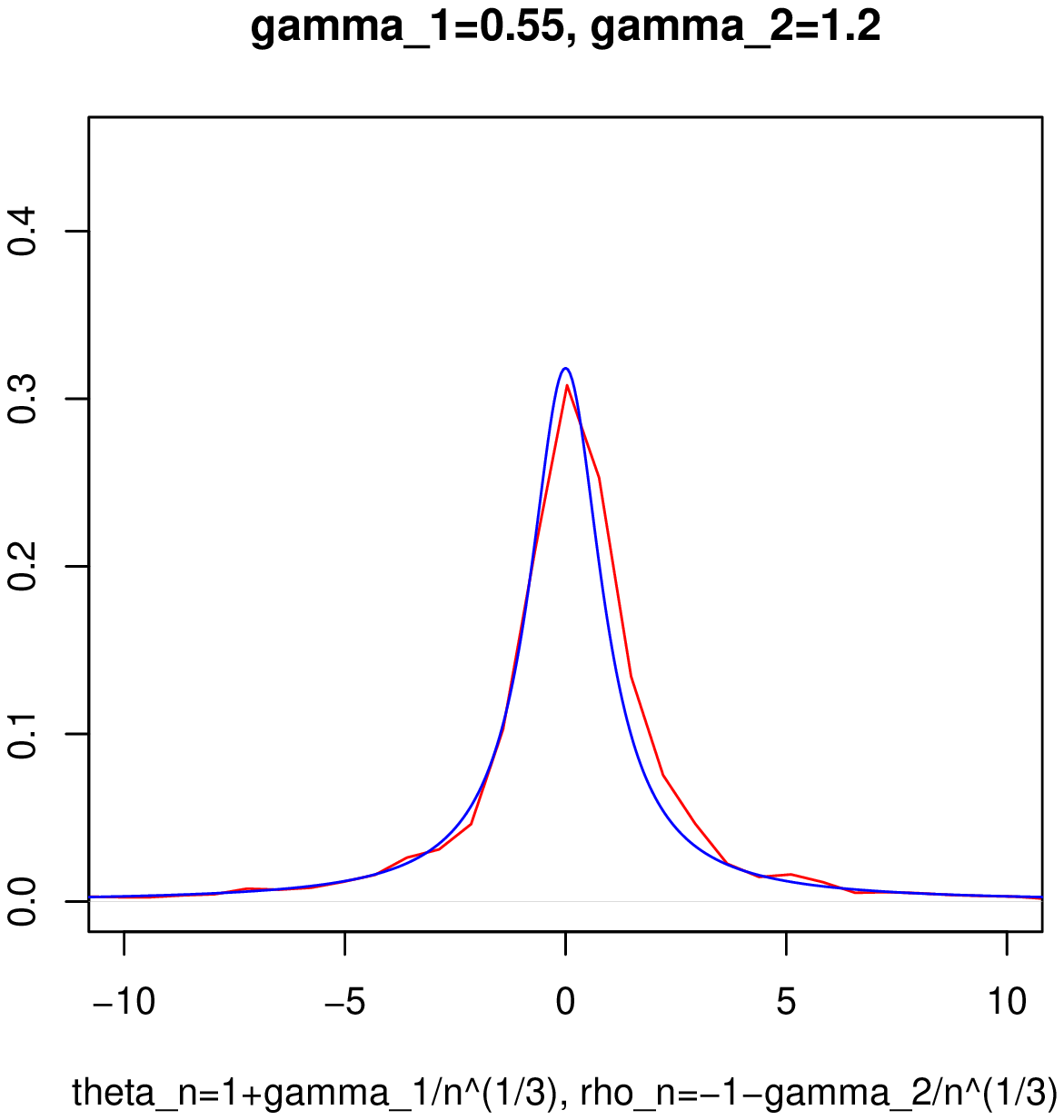}
\includegraphics[width=2in]{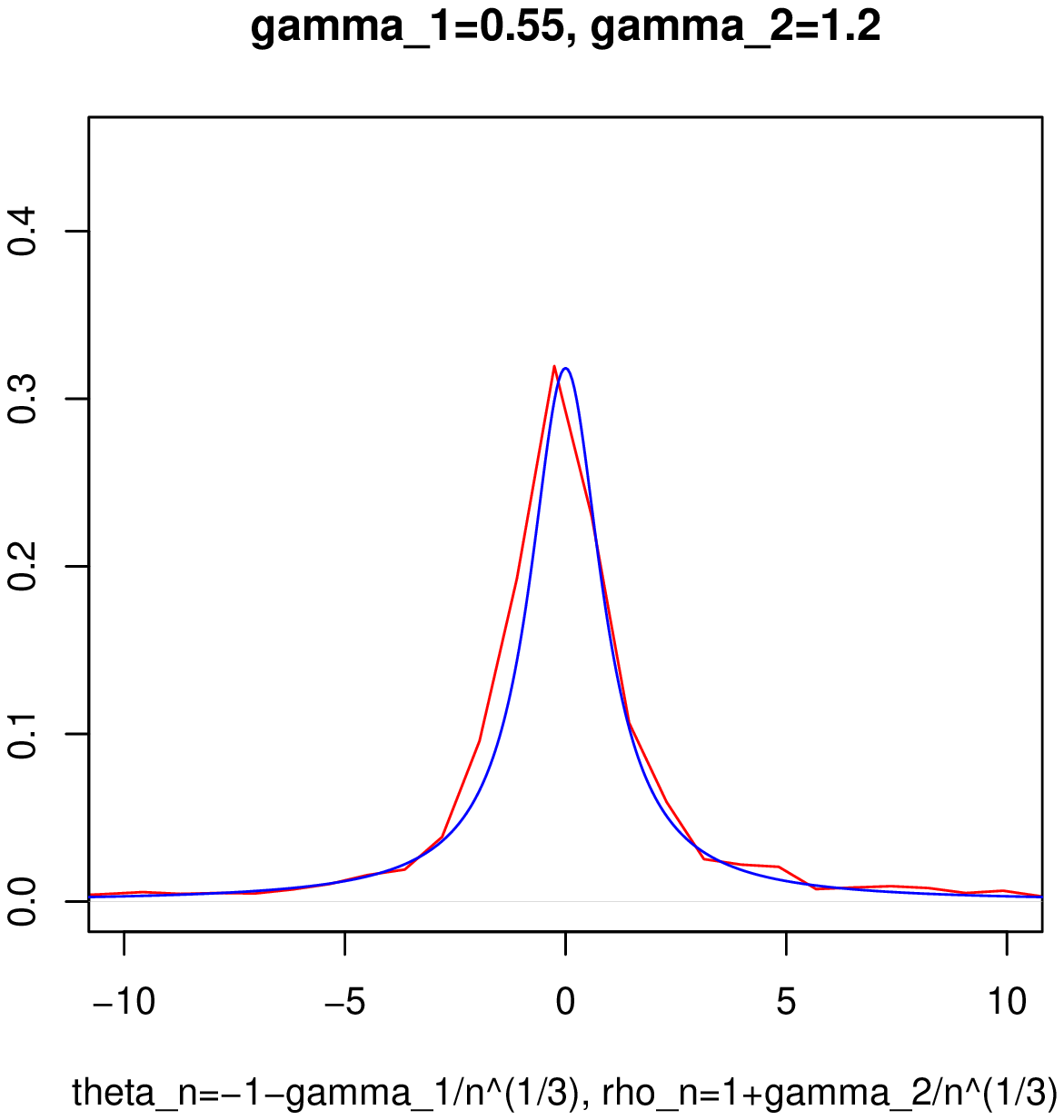}

\section{Proofs of main results}

\subsection{Some preliminary lemmas}

In this subsection, some lemmas are given which play an important role in our following analysis. To obtain the decomposition of $\hat{\theta}_n-\theta_n$, we need introduce some notations. For all $1\leq l\leq n$, let
\begin{align}\label{S-P notion}
P_{l,n}=\sum_{k=1}^l{X_{k,n}X_{k-1,n}}, \quad S_{l,n}=\sum_{k=1}^l{X_{k,n}^2},
\end{align}
and
\begin{align}\label{L-M-N notion}
L_l=\sum_{k=1}^lV_k^2, \quad
M_{l,n}=\sum_{k=1}^l{X_{k-1,n}V_k}, \quad
N_{l,n}=\sum_{k=2}^l{X_{k-2,n}V_k}.
\end{align}
In addition, denote
$$
P_{n,n}:=P_n,\quad S_{n,n}:=S_n, \quad M_{n,n}:=M_n,\quad N_{n,n}:=N_n.
$$
Then we have
\begin{align}\label{decomposition1}
\hat{\theta}_n-\theta_n=\frac{P_n-\theta_nS_{n-1,n}}{S_{n-1,n}}.
\end{align}
Based on the ideas in Bercu \& Pro\"{i}a \cite{Bercu-2013} and Phillips \& Magdalinos \cite{Philllips-Magdalinos}, we deal with the denominator and numerator of (\ref{decomposition1}) respectively. For convenience, define
\begin{equation}\label{def-xieta}
\begin{aligned}
\xi_{\theta_n}&=\frac{1}{\sqrt{k_n}}\sum_{l=1}^n{\theta_n^{-l}V_l},\quad \eta_{\theta_n}=\frac{1}{\sqrt{k_n}}\sum_{l=1}^n{\theta_n^{-(n-l)-1}V_l},\\
\xi_{\rho_n}&=\frac{1}{\sqrt{k_n}}\sum_{l=1}^n{\rho_n^{-l}V_l},\quad \eta_{\rho_n}=\frac{1}{\sqrt{k_n}}\sum_{l=1}^n{\rho_n^{-(n-l)-1}V_l},\\
&\qquad\quad\varphi_{\theta_n}=\frac{1}{n\sqrt{k_n}}\sum_{l=1}^n{(n-l+1)\theta_n^{-l}V_l}.
\end{aligned}
\end{equation}
Then, by some tedious calculations (see Appendix for details), we can write
\begin{align}\label{S_n}
S_{n-1,n}=\frac{1}{\theta_n^2-1}X_{n,n}^2+R_{n1}
\end{align}
and
\begin{align}\label{P_n}
P_n-\theta_nS_{n-1,n}=&\frac{\theta_n\rho_n}{(\theta_n\rho_n-1)(\theta_n-\rho_n)}k_n\theta_n^n\rho_n^n\xi_{\theta_n}\xi_{\rho_n}
\nonumber\\
&+\frac{\rho_n^2}{(\theta_n-\rho_n)(1-\rho_n^2)}k_n\rho_n^{2n}\xi_{\rho_n}^2+R_{n2},
\end{align}
where,
\begin{align}\label{R1}
R_{n1}=&\frac{2\theta_n\rho_n}{(1-\theta_n\rho_n)(\theta_n^2-1)}X_{n,n}\varepsilon_{n,n}
-\frac{\rho_n^2(1+\theta_n\rho_n)}{(1-\theta_n\rho_n)(1-\theta_n^2)(1-\rho_n^2)}\varepsilon_{n,n}^2\nonumber\\
&+\frac{1+\theta_n\rho_n}{(1-\theta_n\rho_n)(1-\theta_n^2)(1-\rho_n^2)}L_n
+\frac{2\theta_n}{(1-\theta_n\rho_n)(1-\theta_n^2)}M_n\nonumber\\
&+\frac{2\rho_n(1+\theta_n\rho_n)}{(1-\theta_n\rho_n)(1-\theta_n^2)(1-\rho_n^2)}\sum_{k=1}^n{\varepsilon_{k-1,n}V_k},
\end{align}
and
\begin{align}\label{R4}
R_{n2}=&\frac{1}{1-\theta_n\rho_n}M_n+\frac{\rho_n}{(1-\theta_n\rho_n)(1-\rho_n^2)}L_n
+\frac{2\rho_n^2}{(1-\theta_n\rho_n)(1-\rho_n^2)}\sum_{k=1}^n{\varepsilon_{k-1,n}V_k}.
\end{align}
Hence, to obtain the asymptotic properties of $S_{n-1,n}$ and $P_n-\theta_nS_{n-1,n}$, we need to deal with the terms appearing in the above equations respectively. And it will be completed by establishing a series of lemmas.

\begin{lem}\label{lem-2}
For model (\ref{model}) with the parameters $|\theta_n|=1+\gamma_1/k_n$ and $|\rho_n|=1+\gamma_2/k_n$, where $\gamma_1,\gamma_2>0$ and $(k_n)_{n\geq1}$ is a sequence of positive numbers increasing to infinity, we have,
\vskip3pt

\noindent(1) if $\theta_n\neq\rho_n$, then

$$
X_{n,n}^2=\frac{\theta_n^2}{(\theta_n-\rho_n)^2}\theta_n^{2n}k_n\xi_{\theta_n}^2
+\frac{\rho_n^2}{(\theta_n-\rho_n)^2}\rho_n^{2n}k_n\xi_{\rho_n}^2
-\frac{2\theta_n\rho_n}{(\theta_n-\rho_n)^2}\theta_n^n\rho_n^nk_n\xi_{\theta_n}\xi_{\rho_n}
$$
and
$$
X_{n,n}\varepsilon_{n,n}=\frac{\theta_n}{\theta_n-\rho_n}\theta_n^n\rho_n^nk_n\xi_{\theta_n}\xi_{\rho_n}
-\frac{\rho_n}{\theta_n-\rho_n}\rho_n^{2n}k_n\xi_{\rho_n}^2,\quad
\varepsilon_{n,n}^2=k_n\rho_n^{2n}\xi_{\rho_n}^2;
$$

\noindent(2) if $\theta_n=\rho_n$, then

$$
X_{n,n}^2=n^2k_n\theta_n^{2n}\varphi_{\theta_n}^2,
\quad X_{n,n}\varepsilon_{n,n}=nk_n\theta_n^{2n}\varphi_{\theta_n}\xi_{\theta_n},
\quad \varepsilon_{n,n}^2=k_n\theta_n^{2n}\xi_{\theta_n}^2.
$$
\end{lem}

\begin{lem}\label{lem-1} Under the aforementioned notations, we have,
\vskip3pt
\noindent(1) if the conditions in Theorem \ref{thm-1} are satisfied, then
$$
\left(\xi_{\theta_n},\eta_{\theta_n},\xi_{\rho_n},\eta_{\rho_n}\right)^{\tau}
\stackrel{\mathcal{L}}\longrightarrow\left(\xi_{\theta}, \eta_{\theta}, \xi_{\rho}, \eta_{\rho}\right)^{\tau},
$$
where $\left(\xi_{\theta}, \eta_{\theta}, \xi_{\rho}, \eta_{\rho}\right)^{\tau}\sim N(0,\Sigma_1)$ with the covariance matrix
\[
\Sigma_1=
\begin{pmatrix}
\frac{\sigma^2}{2\gamma_1} & 0 & \frac{\sigma^2}{\gamma_1+\gamma_2} & 0
\\
0&\frac{\sigma^2}{2\gamma_1}& 0 & \frac{\sigma^2}{\gamma_1+\gamma_2}\\
\frac{\sigma^2}{\gamma_1+\gamma_2} & 0 & \frac{\sigma^2}{2\gamma_2}& 0\\
0 & \frac{\sigma^2}{\gamma_1+\gamma_2} & 0 & \frac{\sigma^2}{2\gamma_2}
\end{pmatrix};
\]
and if the conditions in Theorem \ref{thm-2} are satisfied, then
$$
\left(\xi_{\theta_n},\eta_{\theta_n},\xi_{\rho_n},\eta_{\rho_n}\right)^{\tau}
\stackrel{\mathcal{L}}\longrightarrow\left(\xi_{\theta}, \eta_{\theta}, \xi_{\rho}, \eta_{\rho}\right)^{\tau},
$$
where $\left(\xi_{\theta}, \eta_{\theta}, \xi_{\rho}, \eta_{\rho}\right)^{\tau}\sim N(0,\Sigma_2)$ with the covariance matrix
\[
\Sigma_2=
\begin{pmatrix}
\frac{\sigma^2}{2\gamma_1} & 0 & 0 & 0
\\
0&\frac{\sigma^2}{2\gamma_1}& 0 & 0\\
0 & 0 & \frac{\sigma^2}{2\gamma_2}& 0\\
0 & 0 & 0 & \frac{\sigma^2}{2\gamma_2}
\end{pmatrix};
\]

\noindent(2) if the conditions in Theorems \ref{thm-1} are satisfied, then
$$
\left(\varphi_{\theta_n},\xi_{\theta_n},\eta_{\theta_n}\right)^{\tau}
\stackrel{\mathcal{L}}\longrightarrow
\left(\varphi_{\theta},\xi_{\theta},\eta_{\theta}\right)^{\tau},
$$
where $\left(\varphi_{\theta},\xi_{\theta},\eta_{\theta}\right)^{\tau}\sim N(0,\Gamma)$ with the covariance matrix
$$
\Gamma=
\begin{pmatrix}
\frac{\sigma^2}{2\gamma_1} & \frac{\sigma^2}{2\gamma_1} & 0
\\
\frac{\sigma^2}{2\gamma_1} & \frac{\sigma^2}{2\gamma_1} & 0 \\
0 & 0 & \frac{\sigma^2}{2\gamma_1}
\end{pmatrix},
$$
and
$$
\left(\varphi_{\theta_n}, \frac{n(\xi_{\theta_n}-\varphi_{\theta_n})}{k_n}+\frac{\xi_{\theta_n}}{2\gamma_1}\right)^{\tau}
\stackrel{\mathcal{L}}\longrightarrow
\left(\varphi_{\theta},\zeta_{\theta}\right)^{\tau},
$$
where $\left(\varphi_{\theta},\zeta_{\theta}\right)^{\tau}\sim N(0,\Xi)$ with
the covariance matrix
$$
\Xi=
\begin{pmatrix}
\frac{\sigma^2}{2\gamma_1} & \frac{\sigma^2}{2\gamma_1^2}
\\
\frac{\sigma^2}{2\gamma_1^2}  & \frac{5\sigma^2}{8\gamma_1^3}
\end{pmatrix}.
$$
\end{lem}

\begin{rmk}
Note that $\xi_{\theta_n}$ and $\eta_{\theta_n}$ are always asymptotically independent as well as $\xi_{\rho_n}$ and $\eta_{\rho_n}$, however, $\xi_{\theta_n}$, $\eta_{\theta_n}$, $\xi_{\rho_n}$ and $\eta_{\rho_n}$ are mutually independent when the main regressor and the AR(1) errors in model (\ref{model}) have opposite correlations, i.e. the parameters $\theta_n$ and $\rho_n$ have the opposite signs.
\end{rmk}

\begin{lem}\label{lem-3} Under the aforementioned notations, we have,
\vskip3pt

\noindent(1) if the conditions in Theorem \ref{thm-1} or \ref{thm-2} are satisfied, then
$$
\frac{\sum_{k=1}^n{\varepsilon_{k-1,n}V_k}}{k_n\rho_n^n}
=\xi_{\rho_n}\eta_{\rho_n}+o_p(1),
\quad \frac{\varepsilon_{n,n}^2}{k_n\rho_n^{2n}}=\xi_{\rho_n}^2;
$$

\noindent(2) in the framework of Theorem \ref{thm-1}, if $\gamma_1>\gamma_2>0$, then
$$
\frac{X_{n,n}^2}{k_n^3\theta_n^{2n}}=\frac{1}{(\gamma_1-\gamma_2)^2}\xi_{\theta_n}^2+o_p(1),\quad
\frac{M_n}{k_n^2\theta_n^n}=\frac{1}{\gamma_1-\gamma_2}\xi_{\theta_n}\eta_{\theta_n}+o_p(1)
$$
and
$$
\frac{X_{n,n}\varepsilon_{n,n}}{k_n^2\theta_n^n\rho_n^n}
=\frac{1}{\gamma_1-\gamma_2}\xi_{\theta_n}\xi_{\rho_n}+o_p(1),
$$
and if $\gamma_2>\gamma_1>0$, then
$$
\frac{X_{n,n}^2}{k_n^3\rho_n^{2n}}=\frac{1}{(\gamma_2-\gamma_1)^2}\xi_{\rho_n}^2+o_p(1),\quad
\frac{M_n}{k_n^2\rho_n^n}=\frac{1}{\gamma_2-\gamma_1}\xi_{\rho_n}\eta_{\rho_n}+o_p(1)
$$
and
$$
\frac{X_{n,n}\varepsilon_{n,n}}{k_n^2\rho_n^{2n}}
=\frac{1}{\gamma_2-\gamma_1}\xi_{\rho_n}^2+o_p(1);
$$

\noindent (3) in the framework of Theorem \ref{thm-2}, if $\gamma_1>\gamma_2>0$, then
$$
\frac{X_{n,n}^2}{k_n\theta_n^{2n}}=\frac{1}{4}\xi_{\theta_n}^2+o_p(1),\quad
\frac{M_n}{k_n\theta_n^n}=\frac{1}{2}\xi_{\theta_n}\eta_{\theta_n}+o_p(1)
$$
and
$$
\frac{X_{n,n}\varepsilon_{n,n}}{k_n\theta_n^n\rho_n^n}
=\frac{1}{2}\xi_{\theta_n}\xi_{\rho_n}+o_p(1),
$$
and if $\gamma_2>\gamma_1>0$, then
$$
\frac{X_{n,n}^2}{k_n\rho_n^{2n}}=\frac{1}{4}\xi_{\rho_n}^2+o_p(1),\quad
\frac{M_n}{k_n\rho_n^n}=\frac{1}{2}\xi_{\rho_n}\eta_{\rho_n}+o_p(1)
$$
and
$$
\frac{X_{n,n}\varepsilon_{n,n}}{k_n\rho_n^{2n}}=\frac{1}{2}\xi_{\rho_n}^2+o_p(1);
$$

\noindent (4) in the framework of Theorem \ref{thm-1}, if $\gamma_2=\gamma_1=\gamma>0$, then
$$
\frac{X_{n,n}^2}{n^2k_n\theta_n^{2n}}=\varphi_{\theta_n}^2,
\quad \frac{M_n}{nk_n\theta_n^n}=\varphi_{\theta_n}\eta_{\theta_n}+o_p(1),
\quad \frac{X_{n,n}\varepsilon_{n,n}}{nk_n\theta_n^{2n}}=\varphi_{\theta_n}\xi_{\theta_n}.
$$
\end{lem}

\begin{rmk}
In the part (1) of Lemma \ref{lem-3}, if $\gamma_1=\gamma_2=\gamma>0$, then $\rho_n=\theta_n$, $\xi_{\rho_n}=\xi_{\theta_n}$, and $\eta_{\rho_n}=\eta_{\theta_n}$, so we have
$$
\frac{\sum_{k=1}^n{\varepsilon_{k-1,n}V_k}}{k_n\theta_n^n}=\xi_{\theta_n}\eta_{\theta_n}+o_p(1),
\quad \frac{\varepsilon_{n,n}^2}{k_n\theta_n^{2n}}=\xi_{\theta_n}^2.
$$
\end{rmk}

We are now in a position to provide the asymptotic estimations of $S_{n-1,n}$ and $P_n-\theta_nS_{n-1,n}$ defined as in (\ref{S-P notion}).

\begin{lem}\label{lem-4}
Under the aforementioned notations, we have,
\vskip3pt

\noindent(1) in the framework of Theorem \ref{thm-1}, if $\gamma_1>\gamma_2>0$, then
$$
\frac{S_{n-1,n}}{k_n^4\theta_n^{2n}}
=\frac{1}{2\gamma_1(\gamma_1-\gamma_2)^2}\xi_{\theta_n}^2+o_p(1),
$$
and if $\gamma_2>\gamma_1>0$, then
$$
\frac{S_{n-1,n}}{k_n^4\rho_n^{2n}}
=\frac{1}{2\gamma_2(\gamma_2-\gamma_1)^2}\xi_{\rho_n}^2+o_p(1);
$$

\noindent(2) in the framework of Theorem \ref{thm-2}, if $\gamma_1>\gamma_2>0$, then
$$
\frac{S_{n-1,n}}{k_n^2\theta_n^{2n}}
=\frac{1}{8\gamma_1}\xi_{\theta_n}^2+o_p(1),
$$
and if $\gamma_2>\gamma_1>0$, then
$$
\frac{S_{n-1,n}}{k_n^2\rho_n^{2n}}
=\frac{1}{2\gamma_2(\theta_n-\rho_n)^2}\xi_{\rho_n}^2+o_p(1)
=\frac{1}{8\gamma_2}\xi_{\rho_n}^2+o_p(1);
$$

\noindent(3) in the framework of Theorem \ref{thm-1}, if $\gamma_2=\gamma_1=\gamma>0$, then
$$
\frac{S_{n-1,n}}{n^2k_n^2\theta_n^{2n}}=\frac{1}{2\gamma}\varphi_{\theta_n}^2+o_p(1).
$$
\end{lem}

Before stating the asymptotic estimations of $P_n-\theta_nS_{n-1,n}$, we first give a proposition which shows the direct idea why we consider the asymptotic distributions of $k_n(\hat{\theta}_n-\theta_n)-(\gamma_2-\gamma_1)$, $n(\hat{\theta}-\theta_n)-\theta_n$ and $\hat{\theta}_n-\rho_n$.

\begin{pro}\label{lem-5}
Under the aforementioned notations, we have,
\vskip3pt

\noindent(1) if $\gamma_2>\gamma_1>0$, then in the framework of Theorem \ref{thm-1},
$$
k_n(\hat\theta_n-\theta_n)\stackrel{P}\longrightarrow(\gamma_2-\gamma_1),
$$
and in the framework of Theorem \ref{thm-2},
\[
\hat\theta_n-\theta_n\stackrel{P}\longrightarrow-\theta_n+\rho_n;
\]

\noindent(2) if $\gamma_1=\gamma_2=\gamma>0$, then
$$
n(\hat\theta_n-\theta_n)-\theta_n\stackrel{P}\longrightarrow0.
$$
\end{pro}

Note that, in the following results, because of the bias of the
least squares estimator $\hat{\theta}_n$ when $\gamma_2>\gamma_1>0$, we consider the $P_n-\rho_nS_{n-1,n}$ instead of $P_n-\theta_nS_{n-1,n}$.

\begin{lem}\label{PS}
Under the aforementioned notations, we have,
\vskip3pt

\noindent(1) in the framework of Theorem \ref{thm-1}, if $\gamma_1>\gamma_2>0$, then
\[
\frac{P_n-\theta_nS_{n-1,n}}{k_n^3\theta_n^n\rho_n^n}
=\frac{1}{\gamma_1^2-\gamma_2^2}\xi_{\theta_n}\xi_{\rho_n}+o_p{(1)},
\]
and if $\gamma_2>\gamma_1>0$, then
\begin{equation*}
\frac{P_n-\rho_nS_{n-1,n}}{k_n^3\theta_n^n\rho_n^n}
=\frac{1}{\gamma_2^2-\gamma_1^2}\xi_{\theta_n}\xi_{\rho_n}+o_p{(1)};
\end{equation*}

\noindent(2) in the framework of Theorem \ref{thm-2}, if $\gamma_1>\gamma_2>0$, then
\[
\frac{P_n-\theta_nS_{n-1,n}}{k_n\theta_n^n\rho_n^n}
=\frac{1}{4}\xi_{\theta_n}\xi_{\rho_n}+o_p{(1)},
\]
and if $\gamma_2>\gamma_1>0$, then
\[
\frac{P_n-\rho_nS_{n-1,n}}{k_n\theta_n^n\rho_n^n}
=-\frac{1}{4}\xi_{\theta_n}\xi_{\rho_n}+o_p{(1)}.
\]

\noindent(3) in the framework of Theorem \ref{thm-1}, if $\gamma_2=\gamma_1=\gamma>0$, then
$$
\frac{n(P_n-\theta_nS_{n-1,n})-\theta_nS_{n-1,n}}{nk_n^3\theta_n^{2n}}
=\frac{1}{2\gamma}\varphi_{\theta_n}\cdot
\left(\frac{n}{k_n}(\xi_{\theta_n}-\varphi_{\theta_n})
+\frac{1}{2\gamma}\xi_{\theta_n}\right)
+o_p(1).
$$
\end{lem}

\begin{rmk}
From all the above results, it can be clearly seen that, the orders of the terms in equations (\ref{S_n})-(\ref{R4}) when $\theta_n\rho_n>0$, are higher than that when $\theta_n\rho_n<0$.
\end{rmk}

\subsection{Proof of Theorems \ref{thm-1} and \ref{thm-2}}

Now is the time to give the proofs to our main results.
Because of the similarity, we only prove Theorem \ref{thm-1} in this section.

\vskip5pt

\noindent\textbf{{Proof of part (1) in Theorem \ref{thm-1}}.}
Recall that,
\[
\hat{\theta}_n-\theta_n=\frac{P_n-\theta_nS_{n-1,n}}{S_{n-1,n}},
\]
by the parts (1) of Lemmas \ref{lem-4} and \ref{PS}, we have
\begin{align*}
k_n\theta_n^n\rho_n^{-n}(\hat{\theta}_n-\theta_n)
=\frac{2\gamma_1(\gamma_1-\gamma_2)}{\gamma_1+\gamma_2}\cdot
\frac{\xi_{\rho_n}}{\xi_{\theta_n}}+o_p(1).
\end{align*}
Then the part (1) of Lemma \ref{lem-1}, together with the continuous mapping theorem, yields the part (1) of Theorem \ref{thm-1}. \ \ \ \ \ \ \ \ \ \ \ \ \ \ \ \ \ \ \ \ \ \ \ \ \ \ \ \ \ \ \ \ \ \ \ \ \ \ \ \ \ \ \ \ \ \ \ \ \ \ \ \ \ \ \ \ \ \ \ \ \ \ \ \ \ \ \ \ \ \ \ \  $\Box$

\vskip5pt

\noindent\textbf{{Proof of part (2) in Theorem \ref{thm-1}}.}
According to Proposition \ref{lem-5}, we consider the asymptotic distribution of
$$
k_n(\hat\theta_n-\theta_n)-(\gamma_2-\gamma_1).
$$
By a simple calculation, we can obtain that
\begin{align*}
\hat{\theta}_n-\theta_n-\frac{\gamma_2-\gamma_1}{k_n}
&=\frac{P_n-\rho_nS_{n-1,n}}{S_{n-1,n}}.
\end{align*}
Note that, the parts (1) of Lemmas \ref{lem-4} and \ref{PS} imply that
\[
k_n\rho_n\theta_n^{-n}
\left(\hat{\theta}_n-\theta_n-\frac{\gamma_2-\gamma_1}{k_n}\right)
=\frac{2(\gamma_2-\gamma_1)}{\gamma_1+\gamma_2}\cdot
\frac{\xi_{\theta_n}}{\xi_{\rho_n}}+o_p(1).
\]
Then, using the part (1) of Lemma \ref{lem-1}, and with the aid of the continuous mapping theorem, we complete the proof.
\ \ \ \ \ \ \ \ \ \ \ \ \ \ \ \ \ \ \ \ \ \ \ \ \ \ \ \ \ \ \ \ \ \ \ \ \ \ \ \ \ \ \ \ \ \ \ \ \ \ \ \ \ \ \ \ \ \ \ \ \ \ \ \ \ \ \ \ \ \ \ \ $\Box$

\vskip5pt

\noindent\textbf{{Proof of part (3) in Theorem \ref{thm-1}}.}
According to Proposition \ref{lem-5}, we consider the asymptotic distribution of
$$
n(\hat\theta_n-\theta_n)-\theta_n.
$$
Note that
\[
n(\hat\theta_n-\theta_n)-\theta_n
=\frac{n(P_n-\theta_nS_{n-1,n})-\theta_nS_{n-1,n}}{S_{n-1,n}},
\]
then the parts (3) of Lemmas \ref{lem-4} and \ref{PS} imply that
\[
\frac{n}{k_n}\left(n(\hat\theta_n-\theta_n)-\theta_n\right)
=\frac{\frac{n}{k_n}(\xi_{\theta_n}-\varphi_{\theta_n})
+\frac{1}{2\gamma}\xi_{\theta_n}}{\varphi_{\theta_n}}
+o_p(1).
\]
Therefore, the applications of the part (2) in Lemma \ref{lem-1} and the continuous mapping theorem complete the proof.
\ \ \ \ \ \ \ \ \ \ \ \ \ \ \ \ \ \ \ \ \ \ \ \ \ \ \ \ \ \ \ \ \ \ \ \ \ \ \ \ \ \ \ \ \ \ \ \ \ \ \ \ \ \ \ \ \ \ \ \ \ \ \ \ \ \ \ \ \ \ \ \ $\Box$

%\begin{lem}\label{lem-5}
%As $n\to\infty$, we have the following results.
%
%\noindent(1) If $\gamma_2>\gamma_1>0$,
%$$
%k_n(\hat\theta_n-\theta_n)\stackrel{P}\longrightarrow(\gamma_2-\gamma_1).
%$$%
%
%
%\noindent(2) If $\gamma_1=\gamma_2=\gamma>0$,
%$$
%n(\hat\theta_n-\theta_n)-\theta_n\stackrel{P}\longrightarrow0,
%$$
%where $\stackrel{{P}}\longrightarrow$ denotes the convergence in probability.
%\end{lem}

%\begin{rmk}
% Lemma \ref{lem-5} gives us the direct idea why we consider the asymptotic %distributions of
%  $k_n(\hat\theta_n-\theta_n)-(\gamma_2-\gamma_1)$ %and %$n(\hat\theta_n-\theta_n)-\theta_n$.\quad$\Box$\\
%\end{rmk}

\section{Technical appendix and proofs}

\noindent\textbf{Proof of Lemma \ref{lem-2}.}
For part (1), since for all $1\leq k\leq n$,
\begin{equation}\label{eq-X}
X_{k,n}={\frac{\theta_n}{\theta_n-\rho_n}}\theta_n^k{\sum_{l=1}^k\theta_n^{-l}V_l}-
{\frac{\rho_n}{\theta_n-\rho_n}}\rho_n^k{\sum_{l=1}^k\rho_n^{-l}V_l},
\end{equation}
by a simple calculation, we can write that
\begin{align*}
X_{n,n}^2=&\frac{\theta_n^2}{(\theta_n-\rho_n)^2}\theta_n^{2n}\Big(\sum_{l=1}^n{\theta_n^{-l}V_l}\Big)^2
+\frac{\rho_n^2}{(\theta_n-\rho_n)^2}\rho_n^{2n}\Big(\sum_{l=1}^n{\rho_n^{-l}V_l}\Big)^2\\
&\quad-\frac{2\theta_n\rho_n}{(\theta_n-\rho_n)^2}\theta_n^n\rho_n^n\sum_{l=1}^n{\theta_n^{-l}V_l}\cdot\sum_{l=1}^n{\rho_n^{-l}V_l}\\
=&\frac{\theta_n^2}{(\theta_n-\rho_n)^2}\theta_n^{2n}k_n\xi_{\theta_n}^2+\frac{\rho_n^2}{(\theta_n-\rho_n)^2}\rho_n^{2n}k_n\xi_{\rho_n}^2
-\frac{2\theta_n\rho_n}{(\theta_n-\rho_n)^2}\theta_n^n\rho_n^nk_n\xi_{\theta_n}\xi_{\rho_n}.
\end{align*}
Moreover, using (\ref{eq-X}) and the fact that
\begin{equation}\label{eq-varepsilon}
\varepsilon_{n,n}=\rho_n^{n}\sum_{l=1}^n{\rho_n^{-l}V_l},
\end{equation}
we can obtain
\begin{align*}
X_{n,n}\varepsilon_{n,n}
%&=\frac{\theta_n}{\theta_n-\rho_n}\theta_n^n\rho_n^n\sum_{l=1}^n{\theta_n^{-l}V_l}\cdot\sum_{l=1}^n{\rho_n^{-l}V_l}
%-\frac{\rho_n}{\theta_n-\rho_n}\rho_n^{2n}\sum_{l=1}^n{\rho_n^{-l}V_l}\cdot\sum_{l=1}^n{\rho_n^{-l}V_l}\\
=\frac{\theta_n}{\theta_n-\rho_n}\theta_n^n\rho_n^nk_n\xi_{\theta_n}\xi_{\rho_n}-\frac{\rho_n}{\theta_n-\rho_n}\rho_n^{2n}k_n\xi_{\rho_n}^2.
\end{align*}
Finally, from (\ref{eq-varepsilon}), we know that
$$
\varepsilon_{n,n}^2=\rho_n^{2n}\Big(\sum_{l=1}^n{\rho_n^{-l}V_l}\Big)^2=k_n\rho_n^{2n}\xi_{\rho_n}^2.
$$
which immediately achieves the proof of part (1).

\vskip5pt

For part (2), since $\rho_n=\theta_n$ under the condition $\gamma_1=\gamma_2$, then
\begin{equation}\label{eq-Xn}
X_{k,n}=\theta_n^k\sum_{l=1}^k{(k-l+1)\theta_n^{-l}V_l}.
\end{equation}
Therefore, we can write that
$$
X_{n,n}^2=\theta_n^{2n}\left(\sum_{l=1}^n{(n-l+1)\theta_n^{-l}V_l}\right)^2
=n^2k_n\theta_n^{2n}\varphi_{\theta_n}^2.
$$
Moreover, from (\ref{eq-varepsilon}), (\ref{eq-Xn}) and $\rho_n=\theta_n$, it follows that
$$
\varepsilon_{n,n}^2=\theta_n^{2n}\Big(\sum_{l=1}^n{\theta_n^{-l}V_l}\Big)^2
=k_n\theta_n^{2n}\xi_{\theta_n}^2,
$$
and
$$
X_{n,n}\varepsilon_{n,n}=\theta_n^{2n}\sum_{l=1}^n{(n-l+1)\theta_n^{-l}V_l}\cdot\sum_{l=1}^n{\theta_n^{-l}V_l}
=nk_n\theta_n^{2n}\varphi_{\theta_n}\xi_{\theta_n},
$$
which achieve the proof of part (2). \ \ \ \ \ \ \ \ \ \ \ \ \ \ \ \ \ \ \
\ \ \ \ \ \ \ \ \ \ \ \ \ \ \ \ \ \ \ \ \ \ \ \ \ \ \ \ \ \ \ \ \ \ \ \ \ \
\ \ \ \ \ \ \ \ \ \ \ \ \ \ \ \ \ \ \ \ \ \ \ \ $\Box$

\vskip5pt

\noindent\textbf{Proof of Lemma \ref{lem-1}.}
For part (1), because of the similarity, we only prove the front half part. By the Cram\'{e}r-Wold device (\cite{Kallenberg}, Corollary 5.5), it is sufficient to show that
for any nonzero vector $\upsilon=(\upsilon_1,\upsilon_2,\upsilon_3,\upsilon_4)$,
$$
\upsilon\left(\xi_{\theta_n},\eta_{\theta_n},\xi_{\rho_n},\eta_{\rho_n}\right)^{\tau}\stackrel{\mathcal{L}}\longrightarrow \upsilon\left(\xi_{\theta},\eta_{\theta},\xi_{\rho},\eta_{\rho}\right)^{\tau}.
$$
In fact, we can write that
$$
\upsilon\left(\xi_{\theta_n},\eta_{\theta_n},\xi_{\rho_n},\eta_{\rho_n}\right)^{\tau}=\sum_{l=1}^n{\xi_{nl}},
$$
where
$$
\xi_{nl}:=\frac{1}{\sqrt{k_n}}\left(\upsilon_1\theta_n^{-l}+\upsilon_2\theta_n^{-(n-l)-1}
+\upsilon_3\rho_n^{-l}+\upsilon_4\rho_n^{-(n-l)-1}\right)V_l,\quad 1\leq l\leq n.
$$
Because $\{\xi_{nl}, 1\leq l\leq n\}$ is a sequence of independent and non-identically distributed random variables,
the variance of $\sum_{l=1}^n{\xi_{nl}}$ can be given by
\begin{align*}
E\left(\sum_{l=1}^n{\xi_{nl}}\right)^2
%&=\upsilon %E\Big(\big(\xi_{\theta_n},\eta_{\theta_n},\xi_{\rho_n},\eta_{\rho_n}\big)^{\tau}
%\big(\xi_{\theta_n},\eta_{\theta_n},\xi_{\rho_n},\eta_{\rho_n}\big)\Big)\upsilon^{\tau}\\
=\upsilon\Sigma_n\upsilon^{\tau},
\end{align*}
where
$$
\Sigma_n=\frac{\sigma^2}{k_n}\sum_{l=1}^{n}
\begin{pmatrix}
\theta_n^{-2l} & \theta_n^{-n-1} &(\theta_n\rho_n)^l & \rho_n^{-n-1}(\theta_n/{\rho_n})^{-l}
\\
\theta_n^{-n-1}&\theta_n^{-2n-2+2l}&\theta_n^{-n-1}(\theta_n/{\rho_n})^l &(\theta_n\rho_n)^{-n-1+l}\\
(\theta_n\rho_n)^l & \theta_n^{-n-1}(\theta_n/{\rho_n})^l & \rho_n^{-2l}
& \rho_n^{-n-1}\\
\rho_n^{-n-1}(\theta_n/{\rho_n})^{-l} & (\theta_n\rho_n)^{-n-1+l} & \rho_n^{-n-1}
 & \rho_n^{-2n-2+2l}
\end{pmatrix}.
$$
By simple but tedious calculations, we can obtain that
\begin{equation}\label{eq-var}
E\left(\sum_{l=1}^n{\xi_{nl}}\right)^2=\upsilon\Sigma_n\upsilon^{\tau}
\rightarrow\upsilon\Sigma_1\upsilon^{\tau},
\end{equation}
where the matrix $\Sigma_1$ is defined in Lemma \ref{lem-1}.
Therefore, to prove this lemma, we only need to show the following Lindeberg condition, i.e. for any $\delta>0$
\begin{equation}\label{Lindeberg}
\sum_{l=1}^nE\left(\xi_{nl}^2I_{\{|\xi_{nl}|>\delta\}}\right)\rightarrow0.
\end{equation}
Note that
\begin{align}\label{bound-ex}
\frac{8}{k_n}\sum_{l=1}^n\Big(\upsilon_1^2\theta_n^{-2l}+\upsilon_2^2\theta_n^{-2(n-l)-2}
+\upsilon_3^2\rho_n^{-2l}+&\upsilon_4^2\rho_n^{-2(n-l)-2}\Big)
\nonumber\\
&\longrightarrow \frac{4(\upsilon_1^2+\upsilon_2^2)}{\gamma_1}+\frac{4(\upsilon_3^2+\upsilon_4^2)}{\gamma_2},
\end{align}
 which implies that the left side of (\ref{bound-ex}) is uniformly bounded by a constant $K\in(0,\infty)$. By the following inequality
 $$
 (x+y+z+w)^2\leq8(x^2+y^2+z^2+w^2),
 $$
 the Lindeberg condition can be written as
\begin{align*}
&\sum_{l=1}^n{E(\xi_{nl}^2)I_{\{|\xi_{nl}|>\delta\}}}\\
&\leq\frac{8}{k_n}\sum_{l=1}^n\bigg(\Big(\upsilon_1^2\theta_n^{-2l}+\upsilon_2^2\theta_n^{-2(n-l)-2}
+\upsilon_3^2\rho_n^{-2l}+\upsilon_4^2\rho_n^{-2(n-l)-2}\Big)\\
&\quad\quad\quad
\cdot E\Big(V_l^2I_{\left\{8\Big(\upsilon_1^2\theta_n^{-2l}+\upsilon_2^2\theta_n^{-2(n-l)-2}
+\upsilon_3^2\rho_n^{-2l}+\upsilon_4^2\rho_n^{-2(n-l)-2}\Big)V_l^2>\delta^2k_n\right\}}\Big)\bigg)\\
& \leq K \operatorname*{max}\limits_{1\leq l \leq n} E\left(V_l^2I_{\left\{8\Big(\upsilon_1^2\theta_n^{-2l}+\upsilon_2^2\theta_n^{-2(n-l)-2}
+\upsilon_3^2\rho_n^{-2l}+\upsilon_4^2\rho_n^{-2(n-l)-2}\Big)V_l^2>\delta^2k_n\right\}}\right)\\
&\leq K E\left(V_1^2I_{\left\{V_1^2>\frac{\delta^2k_n}{8(\upsilon_1^2+\upsilon_2^2+\upsilon_3^2+\upsilon_4^2)}\right\}}\right),
\end{align*}
An application of the integrability of $V_1^2$ completes the checking of the Lindeberg condition (\ref{Lindeberg}).

\vskip5pt

Now, we turn to prove part (2). Denote the covariance matrices of $\left(\varphi_{\theta_n},\xi_{\theta_n},\eta_{\theta_n}\right)^{\tau}$
and $\left(\varphi_{\theta_n}, \, \frac{n}{k_n}(\xi_{\theta_n}-\varphi_{\theta_n})
+\frac{1}{2\gamma_1}\varphi_{\theta_n}\right)^{\tau}$ respectively by
$$
\Gamma_n=\frac{\sigma^2}{n^2k_n}\sum_{l=1}^{n}
\begin{pmatrix}
(n-l+1)^2\theta_n^{-2l} & n(n-l+1)\theta_n^{-2l} & n(n-l+1)\theta_n^{-n-1}
\\
n(n-l+1)\theta_n^{-2l}& n^2\theta_n^{-2l} & n^2\theta_n^{-n-1}
\\
n(n-l+1)\theta_n^{-n-1}  & n^2\theta_n^{-n-1} & n^2\theta_n^{-2n-2+2l}
\end{pmatrix}
$$
and
$$
\Xi_n=E\left(\left(\varphi_{\theta_n}, \frac{n(\xi_{\theta_n}-\varphi_{\theta_n})}{k_n}+\frac{1}{2\gamma_1}\xi_{\theta_n}\right)^{\tau}
\left(\varphi_{\theta_n}, \frac{n(\xi_{\theta_n}-\varphi_{\theta_n})}{k_n}+\frac{1}{2\gamma_1}\xi_{\theta_n}\right)\right).
$$
By simple but tedious calculations, we have
\begin{equation*}
\Gamma_n\rightarrow\Gamma,\quad \Xi_n\rightarrow\Xi,
\end{equation*}
where the matrices $\Gamma$ and $\Xi$ are defined in Lemma \ref{lem-1}.
Then, using a similar argument in the proof of (\ref{Lindeberg}), we can
establish the Lindeberg conditions, i.e. for any $\delta>0$, nonzero vectors $\omega=(\omega_1,\omega_2,\omega_3)$ and $\kappa=(\kappa_1,\kappa_2)$,
\begin{equation*}
\sum_{l=1}^nE\left(\psi_{nl}^2I_{\{|\psi_{nl}|>\delta\}}\right)\rightarrow0,
\quad \sum_{l=1}^nE\left(\chi_{nl}^2I_{\{|\chi_{nl}|>\delta\}}\right)\rightarrow0,
\end{equation*}
where for $1\leq l\leq n$,
$$
\psi_{nl}:=\frac{1}{n\sqrt{k_n}}\left(\omega_1(n-l+1)\theta_n^{-l}+\omega_2n\theta_n^{-l}
+\omega_3n\theta_n^{-(n-l)-1}\right)V_l,
$$
and
$$
\chi_{nl}:=\frac{1}{\sqrt{k_n}}\left(\frac{\kappa_1(n-l+1)}{n}
+\kappa_2\left(\frac{l-1}{k_n}+\frac{1}{2\gamma_1}\right)\right)\theta_n^{-l}V_l.
$$
Therefore,  the proof of part (2) can be achieved. \ \ \ \ \ \ \ \ \ \ \ \ \
\ \ \ \ \ \ \ \ \ \ \ \ \ \ \ \ \ \ \ \ \ \ \ \ \ \ \ \ \ $\Box$

\vskip5pt

\noindent\textbf{Proof of Lemma \ref{lem-3}.}
Based on the proof of (10) in Phillips \& Magdalions
\cite{Philllips-Magdalinos}, we can obtain that
$$
\frac{\sum_{k=1}^n{\varepsilon_{k-1,n}V_k}}{k_n\rho_n^n}
=\xi_{\rho_n}\eta_{\rho_n}+o_p(1).
$$
Combined with Lemma \ref{lem-2}, this completes the proof of part (1).

\vskip5pt

We now turn to prove part (2). Because of the similarity of the method, we only deal with the case, $\gamma_1>\gamma_2>0$. First, using Lemma \ref{lem-2}, we have
\begin{align*}
\frac{X_{n,n}^2}{k_n^3\theta_n^{2n}}
%&=\frac{1}{k_n^3\theta_n^{2n}}\left(\frac{\theta_n^2}{(\theta_n-\rho_n)^2}\theta_n^{2n}k_n\xi_{\theta_n}^2
%+\frac{\rho_n^2}{(\theta_n-\rho_n)^2}\rho_n^{2n}k_n\xi_{\rho_n}^2
%-\frac{2\theta_n\rho_n}{(\theta_n-\rho_n)^2}\theta_n^n\rho_n^nk_n\xi_{\theta_n}\xi\rho_n\right)\\
=\frac{\theta_n^2}{(\theta_n-\rho_n)^2k_n^2}\xi_{\theta_n}^2+\frac{\rho_n^{2n+2}}{(\theta_n-\rho_n)^2k_n^2\theta_n^{2n}}\xi_{\rho_n}^2
-\frac{2\theta_n\rho_n^{n+1}}{(\theta_n-\rho_n)^2k_n^2\theta_n^n}\xi_{\theta_n}\xi_{\rho_n}.\\
\end{align*}
Lemma \ref{lem-1}, together with some simple calculations, shows that
\begin{equation}\label{eq-8}
\frac{\theta_n^2}{(\theta_n-\rho_n)^2k_n^2}\xi_{\theta_n}^2=\frac{1}{(\gamma_1-\gamma_2)^2}\xi_{\theta_n}^2+o_p(1).
\end{equation}
Note that ${\rho_n^{n}}{\theta_n^{-n}}=o(1)$.
Combined with (\ref{eq-8}), this gives
\[
\frac{X_{n,n}^2}{k_n^3\theta_n^{2n}}
=\frac{1}{(\gamma_1-\gamma_2)^2}\xi_{\theta_n}^2+o_p(1).
\]
For the estimation of $M_n$, by (\ref{eq-X}), we can obtain that
\begin{align*}
\frac{M_n}{k_n^2\theta_n^n}
&=\frac{1}{k_n^2\theta_n^n}\sum_{k=2}^n{\Big(\frac{\theta_n}{\theta_n-\rho_n}\theta_n^{k-1}\sum_{l=1}^{k-1}{\theta_n^{-l}V_l}
-\frac{\rho_n}{\theta_n-\rho_n}\rho_n^{k-1}\sum_{l=1}^{k-1}{\rho_n^{-l}V_l}\Big)V_k}\\
&=\frac{1}{k_n^2\theta_n^n}\sum_{k=2}^n{\Bigg(\frac{\theta_n}{\theta_n-\rho_n}\theta_n^{k-1}\Big(\sum_{l=1}^n{\theta_n^{-l}V_l}-\sum_{l=k}^n{\theta_n^{-l}V_l}\Big)
V_k\Bigg)}\\
&\quad-\frac{1}{k_n^2\theta_n^n}\sum_{k=2}^n{\Bigg(\frac{\rho_n}{\theta_n-\rho_n}\rho_n^{k-1}\Big(\sum_{l=1}^n{\rho_n^{-l}V_l}-\sum_{l=k}^n{\rho_n^{-l}V_l}\Big)
V_k\Bigg)}\\
&=\frac{\theta_n}{(\theta_n-\rho_n)k_n^2}\sum_{k=2}^n{\theta_n^{-(n-k)-1}V_k}\cdot\sum_{l=1}^{n}{\theta_n^{-l}V_l}\\
&\quad-\frac{\rho_n^n}{\theta_n^n}\cdot\frac{\rho_n}{(\theta_n-\rho_n)k_n^2}\sum_{k=2}^n{\rho_n^{-(n-k)-1}V_k}\cdot\sum_{l=1}^{n}{\rho_n^{-l}V_l}\\
&\quad-\frac{\theta_n}{(\theta_n-\rho_n)k_n^2\theta_n^n}\sum_{k=2}^n{\theta_n^{k-1}V_k}\cdot\sum_{l=k}^{n}{\theta_n^{-l}V_l}\\
&\quad+\frac{\rho_n^n}{\theta_n^n}\cdot
\frac{\theta_n}{(\theta_n-\rho_n)k_n^2\rho_n^n}\sum_{k=2}^n{\rho_n^{k-1}V_k}
\cdot\sum_{l=k}^{n}{\rho_n^{-l}V_l}.
\end{align*}
By a simple calculation, one can see that
\begin{equation}\label{eq-10}
\frac{\theta_n}{(\theta_n-\rho_n)k_n^2}\sum_{k=2}^n{\theta_n^{-(n-k)-1}V_k}\cdot\sum_{l=1}^{n}{\theta_n^{-l}V_l}
=\frac{1}{\gamma_1-\gamma_2}\xi_{\theta_n}\eta_{\theta_n}+o_p(1),
\end{equation}
and
\begin{equation}\label{eq-11}
\frac{\rho_n}{(\theta_n-\rho_n)k_n^2}\sum_{k=2}^n{\rho_n^{-(n-k)-1}V_k}\cdot\sum_{l=1}^{n}{\rho_n^{-l}V_l}
=\frac{1}{\gamma_1-\gamma_2}\xi_{\rho_n}\eta_{\rho_n}+o_p(1).
\end{equation}
Moreover,
\begin{equation}\label{eq-u}
\begin{aligned}
&\frac{\theta_n}{k_n^2(\theta_n-\rho_n)\theta_n^n}\sum_{k=2}^n{\theta_n^{k-1}V_k}\cdot\sum_{l=k}^n{\theta_n^{-l}V_l}\\
&=\frac{\theta_n^{-n}}{k_n^2(\theta_n-\rho_n)}\sum_{k=2}^n{V_k^2}
+\frac{\theta_n^{1-n}}{k_n^2(\theta_n-\rho_n)}\sum_{k=2}^n{\theta_n^{k-1}V_k}\cdot\sum_{l=k+1}^n{\theta_n^{-l}V_l}.
\end{aligned}
\end{equation}
Applying the law of large numbers for the sequence $\{V_k, k\geq1\}$, we obtain that
$$
\frac{\theta_n^{-n}}{k_n^2(\theta_n-\rho_n)}
\sum_{k=2}^n{V_k^2}=O_p\big(\theta_n^{-n}\frac{n}{k_n}\big)=o_p(1).
$$
As for the second term on the right of (\ref{eq-u}),
note that the sequence
$$
\{\sum_{l=k+1}^n{\theta_n^{k-l-1}V_l}V_k,\; 2\leq k\leq n\}
$$
is uncorrelated, which implies that
\begin{align*}
&E\left(\frac{\theta_n}{k_n^2(\theta_n-\rho_n)\theta_n^n}\sum_{k=2}^n{\theta_n^{k-1}V_k}\cdot\sum_{l=k+1}^n{\theta_n^{-l}V_l}\right)^2\\
&=\frac{\sigma^2\theta_n^{2-2n}}{k_n^4(\theta_n-\rho_n)^2}\sum_{k=2}^n{E\bigg(\sum_{l=k+1}^n{\theta_n^{k-l-1}V_l}\bigg)^2}\\
&=\frac{\sigma^4\theta_n^{2-2n}}{k_n^4(\theta_n-\rho_n)^2}\sum_{k=2}^n{\sum_{l=k+1}^n{\theta_n^{2(k-l-1)}}}
%&=\frac{\sigma^4\theta_n^{-2n}}{k_n^4(\theta_n-\rho_n)^2}\cdot\frac{n(1-\theta_n^{-2n})}{\theta_n^2-1}
=O\big(\theta_n^{-2n}\frac{n}{k_n}\big)=o(1).
\end{align*}
Consequently,
\begin{equation}\label{eq-12}
\frac{\theta_n}{k_n^2(\theta_n-\rho_n)\theta_n^n}\sum_{k=2}^n{\theta_n^{k-1}V_k}\cdot\sum_{l=k}^{n}{\theta_n^{-l}V_l}=o_p(1).
\end{equation}
An identical discussion can establish that
\begin{equation}\label{eq-13}
\frac{\theta_n}{k_n^2(\theta_n-\rho_n)\rho_n^n}\sum_{k=2}^n{\rho_n^{k-1}V_k}\cdot\sum_{l=k}^{n}{\rho_n^{-l}V_l}=o_p(1).
\end{equation}
Now, from (\ref{eq-10})-(\ref{eq-13}), it follows that
\begin{align*}
\frac{M_n}{k_n^2\theta_n^n}=\frac{1}{\gamma_1-\gamma_2}\xi_{\theta_n}\eta_{\theta_n}+o_p(1).
\end{align*}

Finally, for $X_{n,n}\varepsilon_{n,n}$, we have no difficulty to obtain by Lemma \ref{lem-2} that
\begin{equation}\label{ex-Xnvarepsilonn}
\begin{aligned}
\frac{X_{n,n}\varepsilon_{n,n}}{k_n^2\theta_n^n\rho_n^n}
&=\frac{\theta_n}{k_n(\theta_n-\rho_n)}\xi_{\theta_n}\xi_{\rho_n}-\frac{\rho_n^{n+1}}{k_n\theta_n^n(\theta_n-\rho_n)}\xi_{\rho_n}^2\\
&=\frac{1}{\gamma_1-\gamma_2}\xi_{\theta_n}\xi_{\rho_n}+o_p(1),
\end{aligned}
\end{equation}
which achieves the proof of part (2).

\vskip5pt

Because the proof of part (3) is similar to that of part (2), we leave it to the interested reader. Finally, we check the part (4). From part (2) of Lemma \ref{lem-2}, it is obvious that
\[
\frac{X_{n,n}^2}{n^2k_n\theta_n^{2n}}=\varphi_{\theta_n}^2
\quad \mbox{and} \quad \frac{X_{n,n}\varepsilon_{n,n}}{nk_n\theta_n^{2n}}=\varphi_{\theta_n}\xi_{\theta_n}.
\]
To estimate $M_n$, firstly, by (\ref{eq-Xn}), we obtain that
\begin{align*}
&\frac{M_n}{nk_n\theta_n^n}=\frac{1}{nk_n\theta_n^n}\sum_{k=2}^n\left(\theta_n^{k-1}\sum_{l=1}^{k-1}{(k-l)\theta_n^{-l}V_l}\right)V_k\\
&=\frac{1}{nk_n}\sum_{k=2}^n{\theta_n^{-(n-k)-1}V_k}\cdot\sum_{l=1}^n{(k-l)\theta_n^{-l}V_l}
-\frac{1}{nk_n\theta_n^n}\sum_{k=2}^n{\theta_n^{k-1}V_k}\cdot\sum_{l=k}^n{(k-l)\theta_n^{-l}V_l}\\
&=\frac{1}{nk_n}\sum_{k=2}^n{\theta_n^{-(n-k)-1}V_k}\cdot\sum_{l=1}^n{\Big((n-l+1)+(k-n-1)\Big)\theta_n^{-l}V_l}\\
&\quad-\frac{1}{nk_n\theta_n^n}\sum_{k=2}^n{\theta_n^{k-1}V_k}\cdot\sum_{l=k}^n{(k-l)\theta_n^{-l}V_l}\\
%&=\frac{1}{nk_n}\sum_{k=2}^n{\theta_n^{-(n-k)-1}V_k}\cdot\sum_{l=1}^n{(n-l+1)\theta_n^{-l}V_l}
%+\frac{1}{nk_n}\sum_{k=2}^n{\theta_n^{-(n-k)-1}V_k}\cdot\sum_{l=1}^n{(k-n-1)\theta_n^{-l}V_l}\\
%&\quad-\frac{1}{nk_n\theta_n^n}\sum_{k=2}^n{\theta_n^{k-1}V_k}\cdot\sum_{l=k}^n{(k-l)\theta_n^{-l}V_l}\\
&=\frac{1}{nk_n}\sum_{k=2}^n{\theta_n^{-(n-k)-1}V_k}\cdot\sum_{l=1}^n{(n-l+1)\theta_n^{-l}V_l}
\\
&\quad+\frac{1}{nk_n}\sum_{k=2}^n{(k-n-1)\theta_n^{-(n-k)-1}V_k}\cdot\sum_{l=1}^n{\theta_n^{-l}V_l}
-\frac{1}{nk_n\theta_n^n}\sum_{k=2}^n{\theta_n^{k-1}V_k}\cdot\sum_{l=k}^n{(k-l)\theta_n^{-l}V_l}.
\end{align*}
By the definitions of $\varphi_{\theta_n}$ and $\eta_{\theta_n}$, we know that
\begin{equation}\label{eq-Mn-1}
\frac{1}{nk_n}\sum_{k=2}^n{\theta_n^{-(n-k)-1}V_k}\cdot\sum_{l=1}^n{(n-l+1)\theta_n^{-l}V_l}
=\varphi_{\theta_n}\eta_{\theta_n}.
\end{equation}
Hence, it is only needed to show that
\begin{equation}\label{eq-Mn-2}
\frac{1}{nk_n}\sum_{k=2}^n{(k-n-1)\theta_n^{-(n-k)-1}V_k}\cdot\sum_{l=1}^n{\theta_n^{-l}V_l}
=o_p(1)
\end{equation}
and
\begin{equation}\label{eq-Mn-3}
\frac{1}{nk_n\theta_n^n}\sum_{k=2}^n{\theta_n^{k-1}V_k}\cdot\sum_{l=k}^n{(k-l)\theta_n^{-l}V_l}
=o_p(1).
\end{equation}
In fact, by the Cauchy-Schwartz inequality and some simple calculations,
\begin{align*}
&E\left|\frac{1}{nk_n}\sum_{k=2}^n{(k-n-1)\theta_n^{-(n-k)-1}V_k}\cdot\sum_{l=1}^n{\theta_n^{-l}V_l}\right|\\
&\leq\frac{1}{nk_n}\left(E\left(\sum_{k=2}^n{(k-n-1)\theta_n^{-(n-k)-1}V_k}\right)^2\right)^{\frac{1}{2}}\cdot
\left(E\left(\sum_{l=1}^n{\theta_n^{-l}V_l}\right)^2\right)^{\frac{1}{2}}\\
&=\frac{\sigma^2}{nk_n}\left(\sum_{i=1}^{n-1}{i^2\theta_n^{-2i}}\right)^{\frac{1}{2}}\cdot
\left(\sum_{l=1}^n{\theta_n^{-2l}}\right)^{\frac{1}{2}}
%&=\frac{\sigma^2}{nk_n}\left(\frac{\theta_n^{-2}}{(1-\theta_n^{-2})^3}
%\left(1+\theta_n^{-2}-n^2\theta_n^{-2(n-1)}+(2n^2-2n-1)\theta_n^{-2n}-(n-1)^2\theta_n^{-2(n+1)}\right)\right)^{\frac{1}{2}}\\
%&\quad\cdot\left(\frac{\theta_n^{-2}-\theta_n^{-2(n+1)}}{1-\theta_n^{-2}}\right)^{\frac{1}{2}}\\
%&=\frac{\sigma^2}{nk_n}\left(O\left(k_n^{3/2}\right)+O\left(\theta_n^{-n}\frac{n}{k_n}\right)\right)
%\cdot O\left(k_n^{1/2}\right)\\
=O\left(\frac{k_n}{n}\right)+O\left(\theta_n^{-n}k_n^{-3/2}\right)=o(1),
\end{align*}
which achieves the proof of (\ref{eq-Mn-2}). As for the checking of (\ref{eq-Mn-3}), since
$$
\frac{1}{nk_n\theta_n^n}\sum_{k=2}^n{\theta_n^{k-1}V_k}\cdot\sum_{l=k}^n{(k-l)\theta_n^{-l}V_l}
=\frac{1}{nk_n\theta_n^n}\sum_{k=2}^n{\theta_n^{k-1}V_k}\cdot\sum_{l=k+1}^n{(k-l)\theta_n^{-l}V_l}
$$
and the sequence $\{\sum_{l=k+1}^n{(k-l)\theta_n^{k-l-1}V_lV_k},~~2\leq k \leq n\}$
is uncorrelated,
we have
\begin{align*}
&E\left(\frac{1}{nk_n\theta_n^n}\sum_{k=2}^n{\theta_n^{k-1}V_k}\cdot\sum_{l=k+1}^n{(k-l)\theta_n^{-l}V_l}\right)^2\\
&=\frac{\sigma^2\theta_n^{-2n}}{n^2k_n^2}\sum_{k=2}^n{E\left(\sum_{l=k+1}^n{(k-l)\theta_n^{k-l-1}V_l}\right)^2}\\
&=\frac{\sigma^4\theta_n^{-2n}}{n^2k_n^2}\sum_{k=2}^n{\sum_{l=k+1}^n{(k-l)^2\theta_n^{2(k-l-1)}}}\\
%&=\frac{\sigma^4\theta_n^{-2n}}{n^2k_n^2}\sum_{k=2}^n{\sum_{i=1}^{n-k}{i^2\theta_n^{-2i-2}}}\\
&\leq\frac{\sigma^4\theta_n^{-2n}}{n^2k_n^2}\sum_{k=2}^n{\sum_{i=1}^{n}{i^2\theta_n^{-2i-2}}}
=O\left(\frac{\theta_n^{-2n}}{nk_n^2}\right)+O\left(\frac{n\theta_n^{-4n}}{k_n^4}\right)=o(1),
\end{align*}
which completes the proof of (\ref{eq-Mn-3}).
Finally, from (\ref{eq-Mn-1})-(\ref{eq-Mn-3}), it follows that
$$
\frac{M_n}{nk_n\theta_n^n}=\varphi_{\theta_n}\eta_{\theta_n}+o_p(1),
$$
which achieves the proof of part (4). \ \ \ \ \ \ \ \ \
\ \ \ \ \ \ \ \ \ \ \ \ \ \ \ \ \ \ \ \ \ \ \ \ \ \ \ \ \ \ \ \ \
\ \ \ \ \ \ \ \ \ \ \ \ \ \ \ \ \ $\Box$

\vskip5pt

\noindent\textbf{Proof of Lemma \ref{lem-4}.}
From (A.14) and (A.23) in Bercu \& Pro\"{i}a \cite{Bercu-2013},
it follows that
\begin{equation}\label{ex-Sn}
\begin{aligned}
&\left(1-(\theta_n+\rho_n)^2-(\theta_n\rho_n)^2\right)S_{n-1,n}\\
&=-X_{n,n}^2-(\theta_n\rho_n)^2X_{n-1,n}^2+L_n
-2\theta_n\rho_n(\theta_n+\rho_n)P_{n-1,n}+2(\theta_n+\rho_n)M_n-2\theta_n\rho_nN_n
\end{aligned}
\end{equation}
and
\begin{equation}\label{ex-Pn}
P_n=\frac{\theta_n+\rho_n}{1+\theta_n\rho_n}S_{n-1,n}+\frac{1}{1+\theta_n\rho_n}M_n+\frac{\theta_n\rho_n}{1+\theta_n\rho_n}X_{n,n}X_{n-1,n,}
\end{equation}
where $P_n$, $S_{n-1,n}$, $L_n$, $M_n$ and $N_n$ are defined as in (\ref{S-P notion}) and (\ref{L-M-N notion}). For some sake of the reader, we list them here again,
\begin{align*}
P_{l,n}=\sum_{k=1}^l{X_{k,n}X_{k-1,n}}, \quad S_{l,n}=\sum_{k=1}^l{X_{k,n}^2},
\end{align*}
and
\begin{align*}
L_l=\sum_{k=1}^lV_k^2, \quad
M_{l,n}=\sum_{k=1}^l{X_{k-1,n}V_k}, \quad
N_{l,n}=\sum_{k=2}^l{X_{k-2,n}V_k}.
\end{align*}
Together with the facts that
$$
N_n=\frac{M_n-\sum_{k=1}^{n}\varepsilon_{k-1,n}V_k}{\theta_n}
$$
and
\begin{equation}\label{ex-XnXn-1}
X_{n-1,n}^{2}=\frac{X_{n,n}^2+\varepsilon_{n,n}^2-2X_{n,n}\varepsilon_{n,n}}{\theta_n},
\quad X_{n-1,n}X_{n,n}=\frac{X_{n,n}^2-X_{n,n}\varepsilon_{n,n}}{\theta_n},
\end{equation}
We can decompose $S_{n-1,n}$ as follows
\begin{equation}\label{eq-15}
S_{n-1,n}=\frac{1}{\theta_n^2-1}X_{n,n}^2+R_{n1},
\end{equation}
where
\begin{align*}
R_{n1}=&\frac{2\theta_n\rho_n}{(1-\theta_n\rho_n)(\theta_n^2-1)}X_{n,n}\varepsilon_{n,n}
-\frac{\rho_n^2(1+\theta_n\rho_n)}{(1-\theta_n\rho_n)(1-\theta_n^2)(1-\rho_n^2)}\varepsilon_{n,n}^2\\
&+\frac{1+\theta_n\rho_n}{(1-\theta_n\rho_n)(1-\theta_n^2)(1-\rho_n^2)}L_n
+\frac{2\theta_n}{(1-\theta_n\rho_n)(1-\theta_n^2)}M_n\\
&+\frac{2\rho_n(1+\theta_n\rho_n)}{(1-\theta_n\rho_n)(1-\theta_n^2)(1-\rho_n^2)}\sum_{k=1}^n{\varepsilon_{k-1,n}V_k}.
\end{align*}
If $\gamma_1>\gamma_2>0$, then by Lemma \ref{lem-3}, we have
$$
\frac{1}{\theta_n^2-1}X_{n,n}^2
%=\frac{1}{\theta_n^2-1}\cdot\frac{k_n^3\theta_n^{2n}}{(\gamma_1-\gamma_2)^2}\xi_{\theta_n}^2+o_p(k_n^4\theta_n^{2n})
=\frac{k_n^4\theta_n^{2n}}{2\gamma_1(\gamma_1-\gamma_2)^2}\xi_{\theta_n}^2+o_p(k_n^4\theta_n^{2n})
$$
and
$$
R_{n1}=o_p(k_n^4\theta_n^{2n}).
$$
Combined with (\ref{eq-15}), this proves the front half part of (1).
If $\gamma_2>\gamma_1>0$, similar to (\ref{eq-15}), we can decompose $S_{n-1,n}$ as follows
\begin{align*}
S_{n-1,n}=&\frac{1}{\theta_n^2-1}X_{n,n}^{2}+\frac{2\theta_n\rho_n}{(1-\theta_n\rho_n)(\theta_n^2-1)}X_{n,n}\varepsilon_{n,n}\\
&-\frac{\rho_n^2(1+\theta_n\rho_n)}{(1-\theta_n\rho_n)(1-\theta_n^2)(1-\rho_n^2)}\varepsilon_{n,n}^2+R_{n3},
\end{align*}
where
\begin{align*}
R_{n3}&=\frac{1+\theta_n\rho_n}{(1-\theta_n\rho_n)(1-\theta_n^2)(1-\rho_n^2)}L_n+\frac{2\theta_n}{(1-\theta_n\rho_n)(1-\theta_n^2)}M_n\\
&\quad+\frac{2\rho_n(1+\theta_n\rho_n)}{(1-\theta_n\rho_n)(1-\theta_n^2)(1-\rho_n^2)}\sum_{k=1}^n{\varepsilon_{k-1,n}V_k}.
\end{align*}
From Lemma \ref{lem-3}, it follows that
\begin{align*}
&\frac{1}{\theta_n^2-1}X_{n,n}^{2}+\frac{2\theta_n\rho_n}{(1-\theta_n\rho_n)(\theta_n^2-1)}X_{n,n}\varepsilon_{n,n}
-\frac{\rho_n^2(1+\theta_n\rho_n)}{(1-\theta_n\rho_n)(1-\theta_n^2)(1-\rho_n^2)}\varepsilon_{n,n}^2\\
&=-\frac{\rho_n^2\rho_n^{2n}k_n}{(\theta_n-\rho_n)^2(1-\rho_n^2)}\xi_{\rho_n}^2+o_p(k_n^4\rho_n^{2n})\\
%\frac{k_n^4\rho_n^{2n}}{2\gamma_1(\gamma_2-\gamma_1)^2}\xi_{\rho_n}^2+\frac{k_n^4\rho_n^{2n}}{\gamma_1(\gamma_2^2-\gamma_1^2)}\xi_{\rho_n}^2
%+\frac{k_n^4\rho_n^{2n}}{2\gamma_1\gamma_2(\gamma_1+\gamma_2)}\xi_{\rho_n}^2+o_p(k_n^4\rho_n^{2n})\\
&=\frac{k_n^4\rho_n^{2n}}{2\gamma_2(\gamma_1-\gamma_2)^2}\xi_{\rho_n}^2+o_p(k_n^4\rho_n^{2n})
\end{align*}
and
$$
R_{n3}=o_p{(k_n^4\rho_n^{2n})},
$$
which achieve the second part of (1).

\vskip5pt

Because the proof of part (2) is similar to that of part (1), we omit it here. Now, only part (3). When $\gamma_1=\gamma_2=\gamma>0$, noting (\ref{eq-15}) and  the fact $\rho_n=\theta_n$, we have
\begin{equation}\label{eq-Sn-1=2}
S_{n-1,n}=\frac{1}{\theta_n^2-1}X_{n,n}^2+R_{n4},
\end{equation}
where
\begin{align*}
R_{n4}=&-\frac{2\theta_n^2}{(1-\theta_n^2)^2}X_{n,n}\varepsilon_{n,n}
-\frac{\theta_n^2(1+\theta_n^2)}{(1-\theta_n^2)^3}\varepsilon_{n,n}^2\\
&+\frac{1+\theta_n^2}{(1-\theta_n^2)^3}L_n
+\frac{2\theta_n}{(1-\theta_n^2)^2}M_n
+\frac{2\theta_n(1+\theta_n^2)}{(1-\theta_n^2)^3}\sum_{k=1}^n{\varepsilon_{k-1,n}V_k}.
\end{align*}
Applying Lemmas \ref{lem-1} and \ref{lem-3}, we obtain that
$$
\frac{1}{\theta_n^2-1}X_{n,n}^2
%=\frac{1}{\theta_n^2-1}n^2k_n\theta_n^{2n}\varphi_{\theta_n}^2
=\frac{n^2k_n^2\theta_n^{2n}}{2\gamma}\varphi_{\theta_n}^2+o_p{(n^2k_n^2\theta_n^{2n})}
$$
and
$$
R_{n4}=o_p(n^2k_n^2\theta_n^{2n}),
$$
which achieve the proof of part (3).
\ \ \ \ \ \ \ \ \ \ \ \ \ \ \ \ \ \ \ \ \ \ \ \
\ \ \ \ \ \ \ \ \ \ \ \ \ \ \ \ \ \ \ \ \ \ \ \ \ \ \ \ \ \
\ \ \ \ \ \ \ \ \ \ \ \ \ \ \ \ \ \ \
$\Box$

\vskip5pt

Because some equations in the proofs of Lemma \ref{PS} will be needed in that of Proposition \ref{lem-5}, we first establish Lemma \ref{PS}.

\vskip5pt

\noindent\textbf{Proof of Lemma \ref{PS}.} Let us begin with the proof of the front half part of (1). By (\ref{ex-Pn}) and (\ref{ex-XnXn-1}), we have
\begin{align*}
P_n-\theta_nS_{n-1,n}
&=\frac{\rho_n(1-\theta_n^2)}{1+\theta_n\rho_n}S_{n-1,n}+\frac{1}{1+\theta_n\rho_n}M_n
+\frac{\theta_n\rho_n}{1+\theta_n\rho_n}X_{n,n}X_{n-1,n}\\
&=\frac{\rho_n(1-\theta_n^2)}{1+\theta_n\rho_n}S_{n-1,n}
+\frac{1}{1+\theta_n\rho_n}M_n+\frac{\rho_n}{1+\theta_n\rho_n}X_{n,n}
\left(X_{n,n}-\varepsilon_{n,n}\right).
\end{align*}
Moreover, from (\ref{eq-15}), it follows that
\begin{equation}\label{ex-PnSn-1}
\begin{aligned}
P_n-\theta_nS_{n-1,n}=&
-\frac{\rho_n}{1-\theta_n\rho_n}X_{n,n}\varepsilon_{n,n}-\frac{\rho_n^3}{(1-\theta_n\rho_n)(1-\rho_n^2)}\varepsilon_{n,n}^2\\
&+\frac{1}{1-\theta_n\rho_n}M_n+\frac{\rho_n}{(1-\theta_n\rho_n)(1-\rho_n^2)}
\left(L_n+2\rho_n\sum_{k=1}^n{\varepsilon_{k-1,n}V_k}\right).
\end{aligned}
\end{equation}
Now, using (1) of Lemma \ref{lem-3} and (\ref{ex-Xnvarepsilonn}), we have
\begin{align*}
&-\frac{\rho_n}{1-\theta_n\rho_n}X_{n,n}\varepsilon_{n,n}-\frac{\rho_n^3}{(1-\theta_n\rho_n)(1-\rho_n^2)}\varepsilon_{n,n}^2\\
&=\frac{\theta_n\rho_n}{(\theta_n\rho_n-1)(\theta_n-\rho_n)}k_n\theta_n^n\rho_n^n\xi_{\theta_n}\xi_{\rho_n}
+\frac{\rho_n^2}{(\theta_n-\rho_n)(1-\rho_n^2)}k_n\rho_n^{2n}\xi_{\rho_n}^2,
\end{align*}
which implies that
\begin{equation}\label{ex-PnSn}
\begin{aligned}
&P_n-\theta_nS_{n-1,n}\\
&=\frac{\theta_n\rho_n}{(\theta_n\rho_n-1)(\theta_n-\rho_n)}k_n\theta_n^n\rho_n^n\xi_{\theta_n}\xi_{\rho_n}
+\frac{\rho_n^2}{(\theta_n-\rho_n)(1-\rho_n^2)}k_n\rho_n^{2n}\xi_{\rho_n}^2+R_{n2},
\end{aligned}
\end{equation}
where
$$
R_{n2}=\frac{1}{1-\theta_n\rho_n}M_n+\frac{\rho_n}{(1-\theta_n\rho_n)(1-\rho_n^2)}L_n
+\frac{2\rho_n^2}{(1-\theta_n\rho_n)(1-\rho_n^2)}\sum_{k=1}^n{\varepsilon_{k-1,n}V_k}.
$$
Since
$$
\frac{\rho_n}{(1-\theta_n\rho_n)(1-\rho_n^2)}L_n=O_p(nk_n^2),
$$
from Lemma \ref{lem-3}, we obtain that
\begin{equation}\label{Rn2}
R_{n2}=o_p{(k_n^3\theta_n^n\rho_n^n\vee k_n^3\rho_n^{2n})}.
\end{equation}
If note the fact,
\begin{equation}\label{eq-23}
\frac{\rho_n^2}{(\theta_n-\rho_n)(1-\rho_n^2)}=\frac{k_n^2}{2\gamma_2(\gamma_2-\gamma_1)}+o(k_n^2),
\end{equation}
then by Lemma \ref{lem-1}, we can get
$$
\frac{\rho_n^2}{(\theta_n-\rho_n)(1-\rho_n^2)}k_n\rho_n^{2n}\xi_{\rho_n}^2=o_p{(k_n^3\theta_n^n\rho_n^n)}.
$$
Finally, from the following fact,
\begin{equation}\label{eq-24}
\frac{\theta_n\rho_n}{(\theta_n\rho_n-1)(\theta_n-\rho_n)}=\frac{k_n^2}{\gamma_1^2-\gamma_2^2}+o(k_n^2),
\end{equation}
we have, by (\ref{Rn2}),
$$
\frac{P_n-\theta_nS_{n-1,n}}{k_n^3\theta_n^n\rho_n^n}=\frac{1}{\gamma_1^2-\gamma_2^2}\xi_{\theta_n}\xi_{\rho_n}+o_p{(1)}.
$$

Now, we turn to the proof of the latter part of (1). Applying the same method as in the proof of (\ref{ex-PnSn}) and combining (\ref{ex-Sn}) and (\ref{ex-Pn}), we can show that
$$
P_n-\rho_nS_{n-1,n}
=\frac{\theta_n\rho_n}{(1-\theta_n\rho_n)(\theta_n-\rho_n)}k_n\theta_n^n\rho_n^n\xi_{\theta_n}\xi_{\rho_n}+R_{n5},
$$
where
\begin{align*}
R_{n5}&=\frac{\theta_n^2}{(\theta_n^2-1)(\theta_n-\rho_n)}k_n\theta_n^{2n}\xi_{\theta_n}^2+\frac{\theta_n}{(1-\theta_n\rho_n)(1-\theta_n^2)}L_n\\
&\quad+\frac{\theta_n^2+1-2\theta_n\rho_n}{(1-\theta_n\rho_n)(1-\theta_n^2)}M_n+\frac{2\theta_n\rho_n}{(1-\theta_n\rho_n)(1-\theta_n^2)}\sum_{k=1}^n{\varepsilon_{k-1,n}V_k}.
\end{align*}
Then, form Lemmas \ref{lem-1} and \ref{lem-3}, we have
$$
R_{n5}=o_p{(k_n^3\theta_n^n\rho_n^n)},
$$
which implies, together with (\ref{eq-24}), that
\begin{equation}\label{eq-25}
\frac{P_n-\rho_nS_{n-1,n}}{k_n^3\theta_n^n\rho_n^n}=\frac{1}{\gamma_2^2-\gamma_1^2}\xi_{\theta_n}\xi_{\rho_n}+o_p{(1)}.
\end{equation}
This completes the proof of f part (1).

\vskip5pt

Although the proof of part (2) is very similar to that of part (1),
there is a little difference on the conditions satisfied by the sequence $(k_n)_{n\geq1}$. Here, we only reproduce the proof of the front half part of (2).
As the same as the proof of the first part of (1), we can get (\ref{ex-PnSn}).
Because $\theta_n$ and $\rho_n$ have opposite signs, we have
$$
\frac{\rho_n}{(1-\theta_n\rho_n)(1-\rho_n^2)}L_n=O_p(nk_n).
$$
Part (3) of Lemma \ref{lem-3} implies that the other two terms in $R_{n2}$ are $O_p(k_n\theta_n^n)$ and $O_p(k_n^2\rho_n^n)$, respectively.
%we can obtain that
%\begin{equation}\label{Rn2+}
%R_{n4}=o_p{(k_n\theta_n^n\rho_n^n\vee k_n^2\rho_n^{2n})}.
%\end{equation}
And Lemma \ref{lem-1} implies that the second term in (\ref{ex-PnSn}) is $O_p{(k_n^2\rho_n^{2n})}$. Note that the following facts,
\begin{equation}
\frac{\rho_n^2}{(\theta_n-\rho_n)(1-\rho_n^2)}=-\frac{k_n}{4\gamma_2}+o(k_n),
\end{equation}
and
\begin{equation}
\frac{\theta_n\rho_n}{(\theta_n\rho_n-1)(\theta_n-\rho_n)}=\frac{1}{4}+o(1),
\end{equation}
then, to obtain the desired result, it is enough to ensure that
\[
n\theta_n^{-n}\rho_n^{-n}\to0 \quad \mbox{and} \quad k_n\rho_n^n\theta_n^{-n}\to0.
\]
Fortunately, $k_n=n^\alpha$ when $\alpha\in(0,1)$ just meets this.

\vskip5pt

For the part (3) of this lemma, i.e. the case $\gamma_1=\gamma_2=\gamma>0$, by (\ref{eq-Sn-1=2}) and (\ref{ex-PnSn-1}), we can write
\begin{align*}
&n\left(P_n-\theta_nS_{n-1,n}\right)-\theta_nS_{n-1,n}\\
&=\frac{n\theta_n}{\theta_n^2-1}X_{n,n}\varepsilon_{n,n}-\frac{\theta_n}{\theta_n^2-1}X_{n,n}^2
+\frac{2\theta_n^3}{(1-\theta_n^2)^2}X_{n,n}\varepsilon_{n,n}-\frac{n\theta_n^3}{(1-\theta_n^2)^2}\varepsilon_{n,n}^2+R_{n6},
\end{align*}
where
\begin{align*}
R_{n6}=&\frac{n}{1-\theta_n^2}M_n+\frac{n\theta_n}{(1-\theta_n^2)^2}\left(L_n+2\theta_n
\sum_{k=1}^n{\varepsilon_{k-1,n}V_k}\right)\\
&+\frac{\theta_n^3(1+\theta_n^2)}{(1-\theta_n^2)^3}\varepsilon_{n,n}^2
-\frac{2\theta_n^2}{(1-\theta_n^2)^2}M_n
-\frac{\theta_n(1+\theta_n^2)}{(1-\theta_n^2)^3}\left(L_n+2\theta_n
\sum_{k=1}^n{\varepsilon_{k-1,n}V_k}\right).
\end{align*}
From Lemma \ref{lem-3}, it follows that
\begin{equation}\label{eq-Rn2}
R_{n6}=O_p\left(k_n^4\theta_n^{2n}\right).
\end{equation}
And using  Lemma \ref{lem-3} again, we can obtain that
\begin{align*}
&\frac{n\theta_n}{\theta_n^2-1}X_{n,n}\varepsilon_{n,n}-\frac{\theta_n}{\theta_n^2-1}X_{n,n}^2
+\frac{2\theta_n^3}{(1-\theta_n^2)^2}X_{n,n}\varepsilon_{n,n}-\frac{n\theta_n^3}{(1-\theta_n^2)^2}\varepsilon_{n,n}^2\\
%&=\frac{\theta_n}{\theta_n^2-1}nk_n\theta_n^{2n}\varphi_{\theta_n}\left(n(\xi_{\theta_n}-\varphi_{\theta_n})
%+\frac{\theta_n^2}{\theta_n^2-1}\xi_{\theta_n}\right)
%+\frac{n\theta_n^3}{(\theta_n^2-1)^2}k_n\theta_n^{2n}\xi_{\theta_n}(\varphi_{\theta_n}-\xi_{\theta_n})\\
&=\frac{\theta_n}{\theta_n^2-1}nk_n\theta_n^{2n}\varphi_{\theta_n}\left(n(\xi_{\theta_n}-\varphi_{\theta_n})
+\frac{k_n}{2\gamma}\xi_{\theta_n}\right)+R_{n7},
\end{align*}
where
$$
R_{n7}=\frac{\theta_n}{\theta_n^2-1}nk_n\theta_n^{2n}\varphi_{\theta_n}\xi_{\theta_n}
\left(\frac{\theta_n^2}{\theta_n^2-1}-\frac{k_n}{2\gamma}\right)
+\frac{n\theta_n^3}{(\theta_n^2-1)^2}k_n\theta_n^{2n}\xi_{\theta_n}(\varphi_{\theta_n}-\xi_{\theta_n}).
$$
However, by part (2) of Lemma \ref{lem-1}, one can see that
\begin{equation}\label{eq-Rn3}
R_{n7}=o_p\left(nk_n^3\theta_n^{2n}\right).
\end{equation}
Above discussions immediately yield the part (3) of this lemma.
  \ \ \ \ \ \ \ \ \ \ \ \ \ \ \ \ \ \ \ \ \ \ \ \ \ \ $\Box$

\vskip5pt

Finally, we end this appendix with the proof of Proposition \ref{lem-5}.

\vskip5pt

\noindent\textbf{Proof of Proposition \ref{lem-5}.} We first deal with the first part of (1) in the proposition. Note that, if $\gamma_2>\gamma_1>0$, then from (\ref{Rn2}), we have
$$
R_{n2}=o_p{(k_n^3\rho_n^{2n})}.
$$
Hence, together with (\ref{eq-23}) and (\ref{eq-24}),
we get
\begin{equation}\label{eq-19}
\frac{P_n-\theta_nS_{n-1,n}}{k_n^3\rho_n^{2n}}=\frac{1}{2\gamma_2(\gamma_2-\gamma_1)}\xi_{\rho_n}^2+o_p{(1)}.
\end{equation}
Consequently, Lemmas \ref{lem-1} and \ref{lem-4} imply that
$$
k_n(\hat\theta_n-\theta_n)\stackrel{P}\longrightarrow(\gamma_2-\gamma_1).
$$

Similarly, we can prove the second part of (1).
As for part (2), if $\gamma_1=\gamma_2=\gamma>0$, from (\ref{ex-PnSn-1}), we can write that
$$
P_n-\theta_nS_{n-1,n}=\frac{\theta_n}{\theta_n^2-1}X_{n,n}\varepsilon_{n,n}+R_{n8},
$$
where
\begin{align*}
R_{n8}=\frac{\theta_n}{(1-\theta_n^2)^2}L_n+\frac{1}{1-\theta_n^2}M_n
+\frac{2\theta_n^2}{(1-\theta_n^2)^2}\sum_{k=1}^n{\varepsilon_{k-1,n}V_k}
-\frac{\theta_n^3}{(1-\theta_n^2)^2}\varepsilon_{n,n}^2.
\end{align*}
Using Lemma \ref{lem-1} and (4) of Lemma \ref{lem-3}, we obtain
$$
R_{n8}=o_p{(nk_n^2\theta_n^{2n})}.
$$
Therefore, we can complete the proof of part (2) in this proposition immediately
 by  Lemma \ref{lem-4} and the fact, $\xi_\theta=\varphi_\theta$. \ \ \ \ \ \ \ \ \ \ \ \ \ \ \ \ \ \ \ \ \ \ \ \
\ \ \ \ \ \ \ \ \ \ \ \ \ \ \ \ \ \ \ \ \ \ \ \ \ \ \ \ \ \
\ \ \ \ \ \ \ \ \ \ \ \ \ \ \ \ \ \ \
$\Box$

\section*{Acknowledgements}

{\small
\noindent The authors wish to express their sincere appreciation to Jianbin Zhao for his kindly and substantive help on the statistical simulations in Section \ref{3}. The work of H. Jiang was partially supported by NSFC (No.11101210), and that of G. Y. Yang was partially supported by NSFC (No. 11201431).
}


\begin{thebibliography}{1234}

\bibitem{Anderson}
Anderson, T. W. (1959). On asymptotic distributions of estimators of parameters of stochastic difference equations. {\it Annals of Mathematical Statistics}, {\bf30}, 676-687.

\bibitem{BasBro1984}
Basawa, I. V. and Brockwell, P. J. (1984). Asymptotic conditional inference for
regular nonergodic models with an application to autoregressive processes. {\it Annals of Statistics}, {\bf12}, 161-171.



\bibitem{Bercu-2013}
Bercu, B. and Pro\"{i}a, F. (2013). A sharp analysis on the asymptotic behavior of thr Durbin-Watson statistic for the First-order autoregressive process.
{\it ESAIM: Probability and Statistics}, {\bf17}, 500-530.


\bibitem{Penda}
Bitseki Penda, V., Djellout, H. and Pro\"{i}a, F. (2013).
Moderate deviations for the Durbin-Watson statistic related to the first-order autoregressive process. {\it ESAIM: Probability and Statistics}, doi: 10.1051/ps/2013038.

\bibitem{Chan2009}
Chan, N. H. (2009). Time series with roots on or near the unit circle. In {\it Iime series: Springer handbooks of financial} (eds T. G. Andersen, R. A. Davis, J. Kreissand, T. Mikosch), New York: Springer. pp. 695-707.

\bibitem{Chan-Wei}
Chan, N. H. and Wei, C. Z. (1987). Asymptotic inference for nearly nonstationary AR(1) processes. {\it Annals of Statistics}, {\bf15}, 1050-1063.

\bibitem{DickFull1979}
Dickey, D. A. and Fuller, W. A. (1979). Distribution of the estimators for autoregressive time series with a unit root. {\it Journal of the American Statistical Association}, {\bf74}, 427-431.

\bibitem{Giraitis}
Giraitis, L. and Phillips, P. C. B. (2006). Uniform limit theory for stationary autoregression. {\it Journal of Time Series Analysis}, {\bf27}, 51-60.

\bibitem{HaszaFuller1979}
Hasza, D. P. and Fuller, W. A. (1979). Estimation for autoregressive processes with unit roots. {\it Annals of Statistics}, {\bf7}, 1106-1120.

\bibitem{JYY}
Jiang, H., Yu, M. M. and Yang, G. Y. (2014). On mildly-stationary second order autoregressive models. {\it Preprint}.

\bibitem{Kallenberg}
Kallenberg, O. (2002). {\it Foundations of Modern Probability}, (2nd edn.), Springer, Berlin.

\bibitem{Magdalinos2012}
Magdalinos, T. (2012). Mildly explosive autoregression under weak and strong dependence. {\it Journal of Econometrics}, 169, 179-187.

\bibitem{Miao-Shen}
Miao, Y. and Shen, S. (2009). Moderate deviation principle for autoregressive processes. {\it Journal of Multivariate Analysis}, 100, 1952-1961.


\bibitem{Miao-Yang}
Miao, Y., Wang, Y. L. and Yang, G. Y. (2014). Moderate deviations principle for empirical covariance from a unit root. {\it Scandinavian Journal of Statistics}, doi: 10.1111/sjos.12104.


\bibitem{Nabeya}
Nabeya, S. and Perron, P. (1994). Local asymptotic distributions related to the AR(1) model with dependent errors. {\it Journal of Econometrics}, 62, 229-264.


\bibitem{Nielsen2009}
Nielsen, B. (2009). Singular vector autoregressions with deterministic terms: strong
consistency and lag order determination. University of Oxford working paper.

%\bibitem{Phil1987}
%Phillips, P. C. B. (1987). Regression theory for near-integrated time series. {\it %Econometrica}, {\bf56}, 1021-1043.

\bibitem{Phillips}
Phillips, P. C. B. (1988). Regression theory for near-integrated time series. {\it Econometrica}, {\bf56}, 1021-1043.


\bibitem{Philllips-Magdalinos}
Phillips, P. C. B. and Magdalinos, T. (2007a). Limit theory for moderate deviations from a unit root. {\it Journal of Econometrics}, 136, 115-130.


\bibitem{Philllips-Magdalinos-1}
Phillips, P. C. B. and Magdalinos, T. (2007b). Limit theory for moderate deviations from a unit root under weak dependence. In: Phillips, G. D. A., Tzavalis, E. (Eds.), {\it The Refinement of Econometric Estimation and Test Procedures}. CUP.

\bibitem{Philllips-Magdalinos-2}
Phillips, P. C. B. and Magdalinos, T. (2013). Inconsistent VAR regression with common explosive roots. {\it Econometric Theory}, {\bf29}, 808-837.

\bibitem{PhiLee2012} Phillips, P. C. B. and Lee, Ji Hyung (2012). VARs with Mixed Roots Near Unity. {\href{http://cowles.econ.yale.edu/P/cd/d18a/d1845.pdf}{\it http://cowles.econ.yale.edu/P/cd/d18a/d1845.pdf}}

\bibitem{Rao1961}
Rao, M. M. (1961). Consistency and limit distributions of estimators of parameters in explosive stochastic difference equations. {\it Annals of Mathematical Statistics}, {\bf32}, 195-218.



\bibitem{Stocker}
Stocker, T. (2007). On the asymptotic bias of OLS in dynamic regression models with autocorrelated errors. {\it Statistical Papers}, 48, 81-93.

\bibitem{White}
White, J. S. (1958). The limiting distribution of the serial correlation coefficient in the explosive case. {\it Annals of Mathematical Statistics}, {\bf29}, 1188-1197.


\end{thebibliography}
\end{document}